\def\ad{{\rm ad}}

\def\End{{\rm End}}

\newtheorem{prop}{Proposition}
\def\M{IM}

\def\R{{\bf R}}

\def\C{{\bf C}}
\def\mod{\hbox{ mod }}

\def\O{{\cal O}}
\def\G{{\cal G}}

\def\ad{{\rm ad}}
\def\tilde{\widetilde}

\def\and{\hbox{\quad and \quad}}

\catcode`\@=12
\newtheorem{theorem}{Theorem}

\documentstyle[12pt]{article}
\input{amssym.def}

\begin{document}

\begin{titlepage}

\rightline{dg-ga/9503016}
\rightline{DAMTP 95-17}
\vskip 0.3cm
\begin{center}
{\LARGE  {\bf Symmetric Monopoles}}
 \vspace*{1cm}
\end{center}

\begin{center}
{\bf {\Large
N.J. Hitchin${}^1$,
N.S. Manton${}^2$, \\}}
\vspace*{.5cm}
{\bf {\Large M.K. Murray${}^3$}}
\end{center}

\vfill

\noindent{\bf Abstract}
We discuss $SU(2)$ Bogomolny monopoles of arbitrary charge $k$
invariant under various symmetry groups. The analysis is
largely in terms of the spectral
curves, the rational maps, and the Nahm equations associated with
monopoles. We consider monopoles
invariant under inversion in a plane, monopoles with cyclic
symmetry, and monopoles having the symmetry of a
regular solid. We introduce the notion of a strongly centred
monopole and show that the space of such monopoles is a geodesic
submanifold of the monopole moduli space.

By solving Nahm's equations we prove the existence of a tetrahedrally
symmetric monopole of charge
$3$ and an octahedrally symmetric monopole of charge $4$, and
determine their spectral curves. Using the geodesic approximation to
analyse the scattering of monopoles with cyclic symmetry, we discover
a novel type of non-planar $k$-monopole scattering process.
\vskip 0.7cm
\centerline{AMS classification scheme numbers: 53C07, 81T13, 32L05}
\vskip 1.3cm
\noindent 1. Department of Pure Mathematics and Mathematical Statistics,
University of Cambridge, 16 Mill Lane, Cambridge CB2 1SB, U.K.
{\em N.Hitchin@pmms.cam.ac.uk}

\noindent 2. Department of Applied Mathematics and Theoretical Physics,
University of Cambridge, Silver Street, Cambridge CB3 9EW, U.K.
{\em N.S.Manton@damtp.cam.ac.uk}

\noindent 3. Department of Pure Mathematics, The University of Adelaide,
Adelaide SA 5005, Australia.  {\em mmurray@maths.adelaide.edu.au}

\end{titlepage}

\bigskip
\section{Introduction}
\label{sec:intro}
In recent years, there has been considerable interest in
monopoles, which are particle-like solitons in a Yang-Mills-Higgs
theory in three spatial dimensions. In this
paper, we shall consider $SU(2)$ Bogomolny monopoles,
which are the finite energy solutions of the Bogomolny
equations (\ref{eq:Bogeq}), \cite{Bog}.  Solutions
are labelled by their magnetic charge, a non-negative integer $k$, and are
physically interpreted as static, non-linear superpositions of $k$
unit charge magnetic monopoles. There is a $4k$-dimensional manifold of
 solutions up to gauge equivalence, known as the $k$-monopole
moduli space $M_k$, and on this there is a naturally defined
Riemannian metric, which is hyperk\"ahler \cite{AtiHit}.

For monopoles moving at modest speeds compared with the speed of
light,  it is a good approximation to model $k$-monopole dynamics by the
geodesic motion on the moduli space $M_k$. This was conjectured some
time ago \cite{Man}, and the consequences explored in some detail
 \cite{AtiHit, GibMan, Woj, BatMon, TemAle}. Very
recently, the validity of the geodesic approximation has been proved
analytically by Stuart \cite{Stu}.

Most studies of Bogomolny monopoles have been concerned either with
the general structure of the $k$-monopole moduli space $M_k$ and its
metric, or with a detailed study of  the geodesics on it for small values
of $k$.
 Little work has
been done on $k$-monopole dynamics for $k>2$. In this paper, we investigate
classes of
$k$-monopole solutions which are invariant under various
symmetry groups and derive results on their scattering. We consider
monopoles invariant under inversion in a
fixed plane, monopoles invariant under a cyclic group of rotations
about a fixed axis, and monopoles invariant
under the symmetry groups of the regular solids, that is, the tetrahedral,
octahedral and icosahedral groups. The existence of $k$-monopoles with cyclic
symmetry was previously shown in \cite{O'RRou}.
Each submanifold of the moduli space $M_k$ consisting of all $k$-monopoles
invariant under a fixed symmetry group is a totally geodesic submanifold.
We therefore obtain various examples of monopole scattering with
symmetry, by finding geodesics on such submanifolds.

Among our most interesting results are proofs of existence of a
tetrahedrally symmetric charge $3$ monopole and an octahedrally
symmetric charge $4$ monopole. We give explicit formulae for the spectral
curves and the
solutions of Nahm's equations corresponding to these monopoles. Our
approach rather strongly indicates that there should be an
icosahedrally symmetric monopole of charge $6$, but a detailed study
of Nahm's equations shows, surprisingly, that such an object does not
exist.

Much of the motivation for the present study of symmetric monopoles came
from results concerning Skyrmions. Skyrmions are $SU(2)$-valued scalar
fields in $\bf R^3$ which minimize Skyrme's energy functional. They have
an integer topological charge $B$, physically identified with baryon
number. Numerical work by Braaten et al. \cite{BraTowCar} has established
that the Skyrmions of charges one to four have, respectively, spherical
symmetry, toroidal symmetry, tetrahedral symmetry and octahedral symmetry.
Similar results have subsequently been obtained using instanton-generated
Skyrmions \cite{AtiMan,LeeMan}. Since a unit charge monopole is
spherically symmetric, and the maximal symmetry of a charge $2$ monopole
is that of a torus, we were led to seek monopoles of charge $3$ and charge
$4$ with tetrahedral and octahedral symmetry, respectively. The Skyrmions
of charge $5$ and charge $6$ are also known, and have rather low symmetry,
so the absence of an icosahedrally symmetric charge $6$ monopole is not so
surprising. We should remark that the relationship between Skyrmions and
monopoles is not systematically understood.  A $B = 1$ Skyrmion has six
degrees of freedom, whereas a unit charge monopole has four.  The moduli
space of charge $k$ monopoles has dimension $4k$.  There is a less
well-defined moduli space of Skyrme fields of baryon number $B$, of
dimension $6B$, and a well-defined space of instanton-generated Skyrme
fields of dimension $8B-1$.  It would be interesting if the charge $B$
monopole moduli space could be identified as a submanifold of either of
these latter spaces.  This is certainly possible for $B = 2$
\cite{AtiMan}.

In Section \ref{sec:monopoles} we review monopoles and their moduli spaces.
Further details of this material can be found in the book \cite{AtiHit}
and in the references contained therein. In Section \ref{sec:ratmaps} we
review the spectral curves and rational maps associated with monopoles. In
Section \ref{sec:invratmap} we show how the rational map changes when a
monopole is inverted in the plane with respect to which the rational map
is defined, and we investigate the monopoles which are invariant under
this inversion. In Section \ref{sec:centratmap} we consider the
holomorphic geometry associated with the centre of a monopole.  We  define
the total phase of a monopole, and introduce the notion of a strongly
centred monopole -- one whose centre is at the origin and whose total
phase is 1. The manifold of strongly centred monopoles is totally geodesic
in $M_k$; in fact up to a $k$-fold covering, it splits off isometrically.

Spectral curves of $k$-monopoles are curves in $TP_1$, the tangent bundle
to the complex projective line, satisfying a number of constraints. In
Section \ref{sec:symmcur} we consider the action of symmetry groups on
general curves in $TP_1$, and present various classes of curves with
cyclic or dihedral symmetry, and with the symmetries of regular solids.
These are candidates for the spectral curves of symmetric monopoles. In
Section \ref{sec:symmellip} we show that the simplest curves with the
symmetries of the regular solids are related to elliptic curves.

In Sections \ref{sec:Nahm} to \ref{sec:icosa} we review the Nahm equations
associated with monopoles and consider the existence of symmetric
monopoles and the corresponding solutions of Nahm's equations. These
equations are in general very difficult to solve explicitly, involving
theta functions of curves of high genus \cite{Erc}. Rather conveniently,
the symmetry conditions imposed here reduce the solutions to ones written
in terms of elliptic functions. We prove thereby the existence of a
tetrahedrally symmetric monopole of charge $3$ (Theorem~\ref{th:tetra}),
and an octahedrally symmetric monopole of charge $4$
(Theorem~\ref{th:octa}), and we determine their spectral curves. We also
prove the {\em non-existence } of an icosahedrally symmetric monopole of
charge $6$ (Theorem~\ref{th:icosa}).

In Section \ref{sec:scatt}, we investigate $k$-monopoles symmetric under
the cyclic group $C_k$. By considering the $C_k$-invariant rational maps,
we show that the strongly centred monopoles with $C_k$ symmetry are
parametrized by a number of geodesic surfaces of revolution in the moduli
space $M_k$. Geodesic motions on these surfaces model either a purely
planar $k$-monopole scattering process, or a novel type of $k$-monopole
scattering, in which $k$ unit charge monopoles collide simultaneously in a
plane and an $l$-monopole and a $(k-l)$-monopole emerge back-to-back along
the line through the $k$-monopole centre, perpendicular to the plane. The
outgoing monopole clusters both become axisymmetric about this line as
their separation increases to infinity. When $k=3$ and $l=1$ or $l=2$,
this geodesic motion passes instantaneously through the tetrahedrally
symmetric 3-monopole (oppositely oriented in the two cases), and when
$k=4$ and $l=2$, through the octahedrally symmetric 4-monopole.

Finally a  warning is necessary for the reader who wishes to delve into
the literature on this subject. There are a number of places in the theory
of monopoles where  one has to make choices and establish conventions.
Most of these are to do with the orientation of $\R^3$ and the induced
complex structure on the twistor space $T P_1$ of all oriented lines in
$\R^3$. Different authors have made different conventions, and minor
inconsistencies can appear to result if the literature is only read in a
cursory manner.

\section{Monopoles}
\label{sec:monopoles}
To define a monopole we start with a pair $(A, \phi)$ consisting of a
connection $1$-form $A$ on ${\bf R}^3$ with values in ${\frak {su}}(2)$, the
Lie algebra of $SU(2)$, and a function $\phi$ (the Higgs field) from
${\bf R}^3$ into ${\frak{su}}(2)$.    The
Yang-Mills-Higgs energy on this pair is
$$
{\cal E}(A, \phi) = \int_{\bf R^3} (|F_A|^2 + |\nabla_A\phi|^2) d^3x
$$
where $F_A= dA  + A\wedge A$ is the curvature of $A$,
$\nabla_A\phi = d\phi + [A, \phi]$ is the
covariant derivative of the Higgs field, and we use
the usual norms on $1$-forms and $2$-forms and the standard
inner product on ${\frak su}(2)$. The
energy is minimized by  the solutions of the Bogomolny equations \cite{Bog}
\begin{equation}
\star F_A = \nabla_A\phi     \label{eq:Bogeq}
\end{equation}
where $\star$ is the Hodge star on forms on ${\bf R}^3$.
These equations, and the energy, are invariant under
gauge transformations, where the gauge group
$\G$ of all maps $g$ from ${\bf R}^3$   to $SU(2)$ acts by
$$
(A, \phi) \mapsto ( gAg^{-1} - dgg^{-1} \, , \, g\phi g^{-1}).
$$
Finiteness of the energy, and the Bogomolny equations, imply certain
boundary conditions at infinity  in ${\bf R}^3$  on the pair $(A, \phi)$
which are spelt out in detail in \cite{AtiHit}. In particular, $|\phi| \to
c$ for some constant $c$ which cannot change with time. Following
\cite{AtiHit}, we fix $c=1$.

A monopole, then, is a gauge equivalence class of
solutions to the Bogomolny equations subject to these boundary conditions.
In some suitable gauge there is a well-defined Higgs  field
at infinity
$$
\phi^\infty \colon S^2_\infty \to S^2 \subset {\frak {su}}(2)
$$
going from the two sphere of all oriented lines through the origin
in $\R^3$ to the unit two-sphere in ${\frak {su}}(2)$.
The degree of $\phi^\infty$ is a positive integer $k$ called
the magnetic charge of the monopole.

Before discussing the moduli space of all solutions
of the Bogomolny equations we need to be a little more precise
and talk about framed monopoles.  We say a pair $(A, \phi)$ is
framed if
$$
\lim_{x_3 \to \infty} \phi(0, 0, x_3) = \pmatrix{ i & 0 \cr 0 & -i
\cr}.
$$
The gauge transformations fixing such pairs are those $g$
with  $\lim_{x_3 \to \infty} g(0, 0, x_3) $ diagonal.
Notice that every monopole can be gauge transformed until it is framed.
So the space of monopoles modulo gauge transformations is the same as
the space of framed monopoles modulo
those gauge transformations that fix them.
We define a framed gauge transformation to be one
such that $\lim_{x_3 \to \infty} g(0, 0, x_3) = 1 $.
The quotient of the set of all framed
monopoles of charge $k$ by the group of
framed gauge transformations is a manifold called the moduli space of
(framed) monopoles of charge $k$ and denoted   $M_k$.  The group of
constant diagonal gauge transformations (a copy of $U(1)$)
acts on $M_k$  and the quotient is called the reduced moduli space $N_k$.
This action is not quite free, because the element $-1$ acts
trivially, but
the group $U(1)/\{\pm 1\}$ acts freely on $M_k$.

The dimension of $M_k$ is $4k$, which can be understood as follows. In the
case that $k=1$ there is a spherically symmetric monopole called the
Bogomolny-Prasad-Sommerfield (BPS) monopole, or unit charge monopole. Its
Higgs field has a single zero at the origin, and its energy density is
peaked there so it is reasonable to think of the origin as the centre or
location of the monopole. The Bogomolny equations are translation
invariant so this monopole can be translated about $\R^3$ and also rotated
by the circle of constant diagonal gauge transformations. This in fact
generates all of $M_1$ which is therefore diffeomorphic to $ S^1\times
\R^3$.   The coordinates on $M_1$ specify the location of the monopole and
what can be thought of as an internal phase. More generally there is an
asymptotic region of the moduli space  consisting of approximate
superpositions of $k$ unit charge monopoles located at $k$ widely
separated points and with $k$ arbitrary phases.

Although it is not possible to assign precisely to a charge $k$ monopole
$k$ points or locations in $\R^3$ it is possible to assign to the monopole
a centre which can be thought of as the average of the locations of the
$k$ particles making up the monopole.   The important property of this
centre is that if we act on the monopole by an isometry of $\R^3$ the
centre moves by the same isometry.   It is also possible to assign to a
$k$-monopole a total phase; this is essentially the product of the phases
of the $k$ unit charge monopoles. If we act on the monopole by a constant
diagonal gauge transformation corresponding to an element $\mu$ of $U(1)$
then the total phase changes by $\mu^{2k}$.

The natural metric on the moduli space $M_k$ is obtained as follows. There
is a flat $L_2$ metric on the space of fields $(A,\phi )$, and this
descends to a curved metric on the space of gauge equivalence classes of
fields. The latter metric, restricted to the $k$-monopole solutions of the
Bogomolny equations is the metric on $M_k$. Since a large part of the
moduli space $M_k$ describes $k$ well-separated unit charge monopoles,
many geodesics on $M_k$ correspond to the scattering of $k$ unit charge
monopoles, and we shall discuss below some particularly symmetric cases of
such scattering.

\section{Spectral curves and rational maps}
\label{sec:ratmaps}
It is not easy to study charge $k$ monopoles directly in terms of their
fields $(A,\phi)$. However, there are various ways of transforming
monopoles to other types of mathematical objects. There is a twistor
theory for monopoles and the result of applying this shows that monopoles
are equivalent to a certain class of holomorphic bundles on the so-called
mini-twistor space $TP_1$. The boundary conditions of the monopole imply
that the holomorphic bundle is determined by an algebraic curve, called
the spectral curve. Monopoles that differ only by a constant diagonal
gauge transformation have the same spectral curve, see \cite{Hit, Hit1}.

The holomorphic bundle of a $k$-monopole is defined  as
follows \cite{Hit}.
Let $\gamma$ be an oriented line in $\R^3$ and let $\nabla_\gamma$
denote covariant differentiation using the connection $A$ along
$\gamma$. One considers the  ordinary differential equation
\begin{equation}
(\nabla_\gamma  - i \phi)v = 0    \label{eq:Hiteq}
\end{equation}
where $v : \gamma \to \C^2$.  The vector space $E_\gamma$ of all solutions
to  equation (\ref{eq:Hiteq}) is two-dimensional and the union of all
these spaces forms a rank two smooth complex vector bundle $E$ over  the
space of all oriented lines in $\R^3$. It can be  shown that this space of
all oriented lines is a complex manifold, in fact isomorphic to $ TP_1$.
One may define on $E$  a holomorphic structure if the monopole satisfies
the Bogomolny equations.  The bundle $E$ has two holomorphic sub-bundles
$E^\pm_1$ whose fibres $(E^\pm_1)_\gamma$ at $\gamma$ are defined to be
the spaces of solutions that decay as $\pm \infty$ is approached along the
line $\gamma$. The set of $\gamma$ where $(E^+_1)_\gamma =
(E^-_1)_\gamma$, so there is a solution decaying at both ends, forms a
curve $S$ in $TP_1$ called the spectral curve of the monopole. It is
possible to show that a decaying solution decays exponentially so the
spectral curve is also the set of all lines along which there is an $L_2$
solution. Intuitively one should think of the spectral lines as being the
lines going through the locations of the monopoles. In the case of charge
$1$, the spectral lines are precisely those going through the centre of
the monopole.

If we describe a typical point in $P_1$ by homogeneous coordinates
$[\zeta_0, \zeta_1]$ then we can cover $P_1$, in the usual way, by two
open sets $U_0$ and $U_1$ where $\zeta_0$ and $\zeta_1$ are non-zero,
respectively. On the set $U_0$ we introduce the coordinate $\zeta =
\zeta_1 / \zeta_0$. Let us also denote by $U_0$ and $U_1$ the pre-images
of these sets under the projection map from $TP_1$ to $P_1$. Then a
tangent vector $\eta {\partial / \partial \zeta}$ at $\zeta$ in $U_0$ can
be given coordinates $(\eta, \zeta)$. These coordinates allow us to
describe an important holomorphic line bundle $L$ on $TP_1$ which has
transition function $\exp(\eta/\zeta)$ on the overlap of $U_0$ and $U_1$.
Similarly for any complex number $\lambda $ we define the bundle
$L^\lambda $ by the transition function $\exp(\lambda \eta/ \zeta)$.
Finally, if $n$ is any integer we define the line bundle $L^\lambda (n)$
to be the tensor product of $L^\lambda $ with the $n$-th power of the
pull-back under projection $TP_1 \to P_1$ of the  dual of the tautological
bundle  on $P_1$. This has transition function $\zeta^{-n}\exp(\lambda
\eta/\zeta)$. The line bundle $L^0$ is clearly trivial so we denote it by
$\O$, and $L^0(n)$ is denoted by $\O(n)$.

Another way of introducing the twistor theory for monopoles is to note, as
in \cite{Mant}, that the Bogomolny equations on ${\bf R}^3$ are equivalent
to the self-dual Yang-Mills equations on ${\bf R}^4$, invariant under
translation in the fourth direction, so monopoles are equivalent to
$S^1$-invariant instantons on $S^1\times {\bf R}^3$. The twistor space $Z$
for $S^1\times {\bf R}^3$ is the quotient of $P_3\backslash P_1$ by a free
${\bf Z}$-action, and is a bundle of groups ${\bf C}^*\times{\bf C}$ over
$P_1$. By the original Atiyah-Ward construction, an instanton corresponds
to a holomorphic bundle on the twistor space $Z$, and if it is
$S^1$-invariant, it descends to $Z/{\bf C}^*=TP_1$. The space $Z$ itself
is the total space of the principal bundle for the line bundle $L$ defined
above by transition functions.

To avoid the potential ambiguity in what we mean by `transition function'
let us be more explicit. The line bundle $L^\lambda (n)$ has non-vanishing
holomorphic sections $\chi_0$ and $\chi_1$ over $U_0$ and $U_1$
respectively and for points in $U_0 \cap U_1$ these satisfy
\begin{equation}
\chi_0 = \zeta^{-n} \exp({\lambda \eta\over \zeta}) \chi_1.  \label{eq:trans}
\end{equation}
If we consider  an arbitrary holomorphic section $f$ of this line bundle
its restriction to $U_0$ and $U_1$ can be written as
$f = f_0 \chi_0 $ and $f= f_1 \chi_1$ respectively
 where $f_0$ and
$f_1$ are holomorphic functions on $U_0$ and
$U_1$.
As a consequence of equation (\ref{eq:trans}) these functions must satisfy
\begin{equation}
f_0 = \zeta^{n} \exp({-\lambda \eta\over \zeta}) f_1   \label{eq:sec}
\end{equation}
at points in the  intersection $U_0 \cap U_1$.

With these definitions we can present the results
that we need.  The sub-bundles
$E_1^\pm$ satisfy $E_1^\pm \simeq L^{\pm 1}(-k)$
and the quotients satisfy $E/E_1^\pm \simeq L^{\mp 1}(k)$.
For a framed monopole there are explicit isomorphisms
so we shall write $=$ instead of $\simeq$. The curve $S$
is defined by the vanishing of the map $ E_1^+ \to E/E_1^-$ and hence
by a section of $(E_1^+)^* \otimes E/E_1^- = \O(2k)$.
In terms of the coordinates $(\eta, \zeta)$, $S$ is defined by an
equation of the form
\begin{equation}
P(\eta, \zeta) \equiv \eta^k + \eta^{k-1} a_1(\zeta)  + \dots + \eta
a_{k-1}(\zeta)
+ a_k(\zeta)= 0,     \label{eq:spec}
\end{equation}
where, for $1 \leq r \leq k$, $a_r(\zeta)$ is a polynomial in $\zeta$
of degree at most $2r$.

The space $TP_1$ has a real structure $\tau$, namely, the anti-holomorphic
involution defined by reversing the orientation of the lines in
$\R^3$. In coordinates it takes the form $\tau(\eta, \zeta) =
(-\bar\eta/\bar\zeta^2,
-1/\bar\zeta)$.  The curve $S$ is fixed by this involution, so we say
that it is real. The reality of $S$ implies that for $1 \leq r \leq
k$,
\begin{equation}
a_r(\zeta) = (-1)^r \zeta^{2r} \overline{a_r(-{1 \over
{\bar\zeta}})}.  \label{eq:real}
\end{equation}
If $k=1$ the spectral curve has the form
$$
\eta = (x_1 + i x_2) -2x_3\zeta - (x_1 - i x_2)\zeta^2
$$
where $x = (x_1, x_2, x_3)$ is any point in $\R^3$, [\cite{Hit},
eq. (3.2)]. Such a curve is called a
real section as it defines a section of the bundle $TP_1 \to P_1$, and
is real in the sense given above.
In terms of the geometry of $\R^3$ this curve is the set of all
oriented lines through the point $x$, so it is the spectral curve of a BPS
monopole located at $x$. We refer to this curve as the ``star'' at $x$.

In \cite{Hit, Hit1} one can find listed all the properties that a curve in
$TP_1$ has
to satisfy to be a spectral curve.  We are interested in
one of these here. From the definition of the spectral curve
we see that over the spectral curve
the line bundles $E_1^+$ and $E_1^-$ coincide as sub-bundles of $E$; in
particular they must be isomorphic. This is equivalent to
saying that  the line bundle
$E_1^+ \otimes (E_1^-)^* = L^2$ is trivial over the curve
or that it admits a non-vanishing
holomorphic section $s$. The real structure $\tau$ can
be lifted to an anti-holomorphic, conjugate linear
map between the line bundles $L^2$ and $L^{-2}$ and hence
the section $s$ can be conjugated to define a new (holomorphic) section
$\tau(s) = \tau \circ s \circ \tau$ of $L^{-2}$ over
$S$. Tensoring these defines a section
$\tau(s)s$ of $L^{-2} \otimes L^2 = \O$
and because $S$ is compact and connected this is a constant.
Because of the framing this constant will be $1$.
Notice that given only $S$ and the fact that $L^2$ is trivial
over $S$, if we can choose a section $s$ such that $\tau(s)s=1$
 then it is unique up to multiplication by a scalar of
modulus one. This circle ambiguity in the choice of
$s$ corresponds to the framing of the monopole.
In fact, let $\mu$ be a complex number of modulus one corresponding
to a constant diagonal gauge transformation with diagonal
entries $\mu$ and $\mu^{-1}$.
Then it is possible to follow through the proof in (\cite{Hit}, pp.
593-4)
and show that if we phase rotate a framed  monopole by $\mu$,
the isomorphism $E^+_1 \to L(-k)$ is multiplied by $\mu$
and  the isomorphism
$E^-_1 \to L^*(-k)$
is multiplied  by $\mu^{-1}$. The section $s$ of
$E_1^+ \otimes (E_1^-)^* = L^2$ is therefore  multiplied by
$\mu^2$. Notice that this is consistent with the fact that
the group $U(1)/\{\pm 1\}$  acts freely on the moduli space $M_k$
of framed monopoles.

The rational map of a monopole was originally described by
Donaldson in terms of solutions to Nahm's equations \cite{Don}.
Hurtubise then showed how it relates to scattering in $\R^3$ and to the
spectral curve of the monopole \cite{Hur}.
It will be convenient for our purposes
to use the  description in terms of spectral curves.

The rational map of a charge $k$ monopole is from $\C$ to $\C \bigcup
\infty$, and
is simply a polynomial $p$ of degree
less than $k$ divided by a monic (leading coefficient = 1) polynomial $q$
of degree $k$  which has no factor
in common with $p$,
$$
R(z) = {p(z) \over q(z)} .
$$
We shall denote by  $R_k$ the space of all these based rational maps.
Donaldson has proved that any such rational
map arises from some unique charge $k$ monopole \cite{Don}, so $R_k$
is diffeomorphic to $M_k$.
The disadvantage of characterising a monopole by its rational map is that
the definition of the map
requires choosing a line and an orthogonal plane in $\R^3$, and this breaks
the symmetries of the problem.  Whereas the Bogomolny
equations are invariant under all the isometries of $\R^3$,
the transformation to a rational map commutes only with
those isometries that preserve the direction of the line.

To define the rational map we fix the fibre $F$ of $TP_1 \to P_1$ where
$\zeta = 0$ and identify it with $\C$. The fibre consists of all lines in
the $x_3$-direction. This corresponds to picking an orthogonal splitting
of $\R^3$ as $\C \times \R$. Each point $z$ in $\C$ is identified with a
point in $F$ by setting $z = \eta$, and hence with an oriented line, the
line $\{(x_1,x_2,x_3) \mid x_3 \in \R\}$ with $z = x_1 + ix_2$. The
intersection of $F$ with $S$ defines $k$ points counted with multiplicity
and  $q(z)$ is defined to be the unique monic polynomial of degree $k$
which has these $k$ points as its roots. Thus $q(z) = P(z,0)$, where $P$
is given by eq.~(\ref{eq:spec}). Recall from (\ref{eq:sec}) that a
holomorphic section $s$ of the bundle $L^2$ is determined locally by
functions $s_0$ and $s_1$, on $U_0 \cap S$ and $U_1 \cap S$ respectively,
such that
\begin{equation}
s_0(\eta, \zeta) =
\exp({-2\eta \over \zeta}) s_1(\eta, \zeta). \label{eq:secL2}
\end{equation}
Let $p(z)$ be the unique polynomial of degree $k-1$ such that
$p(z)  = s_0(z , 0) \mod q(z)$.
The rational map of the monopole is then $R(z) = p(z)/q(z)$.
If the roots of $q(z)$ are distinct complex numbers
$\beta_1, \dots , \beta_k$
then the polynomial $p(z)$ is determined by its values
$p(\beta_i) = s_0(\beta_i, 0)$ for all $i = 1 , \dots , k$.

Let us make a brief remark about the construction of the rational
map as scattering data  in $\R^3$. More details are given
in \cite{Hur} and  \cite{AtiHit}.
The points where $S$ intersects $F$, the zeros of $q$, correspond to the
lines in the $x_3$-direction admitting a solution of eq.(\ref{eq:Hiteq})
decaying at both ends.
Assume these lines are distinct and label them by the corresponding
complex numbers $\beta_i$. Pick for each line a solution $v(\beta_i, x_3)$
decaying at both ends. In the
regions where $x_3$ is large positive and  large negative
there are choices of asymptotically flat gauge such that
$$
 \lim_{x_3 \to \infty} (x_3)^{-k/2} e^{x_3} v(\beta_i, x_3)  =
v_i^+ \pmatrix{1\cr 0\cr}
$$
and
$$
 \lim_{x_3 \to - \infty} (x_3)^{-k/2} e^{-x_3} v(\beta_i, x_3)  =
 v_i^- \pmatrix{0\cr1\cr}.
$$
The rational
map with our conventions is determined by
$$
p(\beta_i) = {v^+_i \over v^-_i}.
$$
This agrees with the results stated in Chapter 16 of ref.\cite{AtiHit},
although Hurtubise's conventions give $p(\beta_i) = v^-_i / v^+_i$.

We deduce from these formulae the action of certain isometries on the
rational maps of monopoles. Let $\lambda \in U(1)$  and  $w \in \C$ define
a rotation and translation, respectively, in the plane $\C$. Let $t \in
\R$ define a translation perpendicular to the plane and let  $\mu \in
U(1)$ define a constant diagonal gauge transformation.  A rational map
$R(z)$ then transforms under the composition of all these transformations
to
$$
\tilde R(z) = \mu^2 \exp(2t) \lambda^{-2k} R(\lambda^{-1}(z-w)).
$$
Note that this is slightly different to the action described
in [\cite{AtiHit}, eq. (2.11)], because of different conventions.

\section{Inverting monopoles}
\label{sec:invratmap}
Consider the inversion  map $I \colon \R^3 \to \R^3$ defined by
$I(x_1, x_2, x_3) = (x_1, x_2, -x_3)$.
This inverts $\R^3$ in the $(x_1, x_2)$ plane.  The inversion map
reverses orientation, and so induces an {\it anti-holomorphic}  map on
the twistor space $TP_1$  which we shall denote by the
same symbol and which in the standard coordinates on $TP_1$ is
$$
I(\eta, \zeta) = ({-\bar\eta \over \bar\zeta^2}, {1\over \bar\zeta}).
$$
To see this note that the real section defined by the
point $I(x_1, x_2, x_3)$ has equation
$$
\eta = (x_1 + i x_2) + 2 x_3 \zeta - (x_1 - ix_2)\zeta^2.
$$
So a point $I(\eta, \zeta) $ is on this curve
if and only if
$$
{-\bar\eta \over \bar\zeta^2} = (x_1 + ix_2)  + 2 x_3 {1\over \bar\zeta} -
   (x_1 - ix_2){1\over \bar\zeta^2}.
$$
Conjugating this equation and clearing the denominators we recover
$$
\eta = (x_1 + ix_2)  - 2 x_3 \zeta - (x_1 - ix_2)\zeta^2,
$$
the equation of the real section defined by the point $(x_1, x_2, x_3)$.
This confirms the formula for $I$.
Notice that $I$ is very similar to the real structure $\tau$; in fact
$I\circ \tau (\eta, \zeta) = (\eta, -\zeta).$

If we invert the monopole defined by the spectral curve $S$ and
section $s$ we obtain
a new curve $I(S)$  and a new section $I(s)$. The definition of
$I(S)$ is straightforward; it is just the image of $S$ under the map
$I$. We shall consider $I(s)$ in a moment. Because $\tau(S) = S$
it follows that $(\eta, \zeta)\in I(S)$ precisely when
$(\eta, -\zeta) \in S$. In particular, the intersection
of $I(S)$ and the fibre $F$ over $\zeta = 0$ is just the intersection of $S$
and $F$.
So if we denote by $I(p)$ and $I(q)$ the numerator and denominator of
the rational map
for the inverted monopole, we see that $I(q) = q$.

Now consider the section $s$. Notice that both $\tau$ and $I$ interchange
the two coordinate patches $U_0$ and $U_1$. The section $\tau(s)$ is
defined locally by
\begin{equation}
\tau(s)_0(\eta, \zeta) = \bar{s}_1(\tau(\eta, \zeta))
\quad , \quad   \tau(s)_1(\eta, \zeta) =
\bar{s}_0(\tau(\eta, \zeta))  \label{eq:tau}
\end{equation}
and the section $I(s)$ by
\begin{equation}
I(s)_0(\eta, \zeta) = \bar{s}_1(I(\eta, \zeta))
\quad , \quad   I(s)_1(\eta, \zeta) =
\bar{s}_0(I(\eta, \zeta)).  \label{eq:I}
\end{equation}
Hence $I(p)$ is defined by
$$
I(p)(z) = I(s)_0(z , 0) \mod q(z) = \bar{s}_1\circ \tau(z , 0) \mod q(z)
$$
using the fact that $\tau(\eta, 0) = I(\eta, 0)$. From
the relation $\tau(s)s = 1$ and (\ref{eq:tau}) it follows that
$(\bar{s}_1\circ \tau) s_0 = 1$ and hence
\begin{equation}
(I(p)p)(z)  =  (\bar{s}_1\circ \tau(z , 0)) s_0(z ,0) \mod q(z)
       =  1 \mod q(z).
\label{eq:Irat}
\end{equation}
Eq.~(\ref{eq:Irat}), and the requirement that the degree of $I(p)$ is less
than $k$, determine $I(p)$ uniquely. If the roots of $q$ are the distinct
complex numbers $\beta_1, \dots, \beta_k$, a useful alternative way of
obtaining $I(p)$ is to notice that it is  the unique polynomial of degree
less than $k$ such that $I(p)(\beta_i) p(\beta_i) = 1$ for all $i = 1,
\dots, k$.

We summarize these results as:
\begin{prop} Let a monopole have spectral curve $P(\eta , \zeta) = 0$
and rational map $p/q$. The inverted monopole has spectral curve
$P(\eta , -\zeta) = 0$ and rational map $I(p)/q$, where $I(p)p = 1 \mod q$.
\end{prop}

It is interesting to consider the subset of monopoles that are invariant
under inversion. Their spectral curves are given by polynomials $P(\eta ,
\zeta)$ which are even in $\zeta$. Their rational maps satisfy $p^2 = 1$
mod $q$, so that $I(p) = p$. This fixed-point set is described by the
following: \begin{prop} The moduli space $IM_k$ of $k$-monopoles invariant
under inversion is a totally geodesic submanifold of $M_k$ of dimension
$2k$. It has $(k+1)$ connected components $IM^m_k$ for $0\le m \le k$. The
component $IM^m_k$ is diffeomorphic to the set of coprime pairs $(r,s)$ of
monic polynomials of degree $m$ and $(k-m)$ respectively.
\end{prop}
\noindent {\it Proof}:
It is a standard fact from differential geometry that the fixed point set
of a  finite group of isometries of a Riemannian manifold is a totally
geodesic submanifold.  For any rational function $R=p/q$ defining a
$k$-monopole, consider $$f(R)=\sum_ip(\beta_i)$$ This is a symmetric
polynomial in the roots of $q$ and hence is polynomial in the coefficients
of $p$ and $q$, and so is continuous on $M_k$. On $IM_k$, $p^2=1 \mod q$,
hence $p(\beta_i)=\pm1$. Thus, restricted to $IM_k$, $f$ takes the values
$k,k-2,\dots,-k$, and we define $IM^m_k=f^{-1}(k-2m)$.

If $p^2-1=0\mod q$ then $q$ divides $(p-1)(p+1)$. Since these factors are
coprime, any irreducible factor of $q$ divides one or the other. Hence we
have monic polynomials $r,s$ of degrees $m,k-m$ respectively, such that
$q=rs$, and $p+1=2ar$, $p-1=2bs$ for polynomials $a,b$. Hence $p=ar+bs$
and $ar-bs=1$.

Conversely, given two coprime monic polynomials $r,s$, the division
algorithm implies that there exist polynomials $a,b$ such that $ar-bs=1$,
and moreover $a$ can be chosen uniquely to have degree less than that of
$s$. Now define $p=ar+bs$. This has degree less than $k$ and
\begin{equation}
p^2=(ar+bs)^2=1+4abq
\label{pq}
\end{equation}
where $q=rs$ is monic of degree $k$. Clearly, from (\ref{pq}), $p$ and $q$
are coprime and so define a rational map $R$.

Now the space of pairs $(r,s)$ of coprime polynomials is the complement of
a hypersurface in ${\bf C}^m\times{\bf C}^{k-m}$ and so is a connected
$2k$-dimensional manifold. Moreover, when $r,s$ both have distinct roots,
the roots of $q$ are $\beta_1,\dots,\beta_m$  and
$\beta_{m+1},\dots,\beta_k$ (the roots of $r$ and $s$ respectively). Hence
on this manifold $f(R)=k-2m$, and so $IM^m_k$ is connected, and as
described in the  proposition.
 \vskip .5cm
 Note that $\M_k^m$ and
$\M_k^{k-m}$ are isomorphic; one is obtained from the
other by multiplying $p$ by $-1$. The simplest of the components is
$\M_k^0$. Here $p(z) \equiv 1$, so the rational maps are of the
form
$$
R(z) = {1 \over q(z)}.
$$
The space $\M_k^0$ is naturally diffeomorphic to the moduli space of $k$
flux vortices in the critically coupled abelian Higgs model, since
$k$-vortex solutions are also parametrised by a single monic polynomial of
degree $k$ \cite{Tau}. However, the metrics in the monopole and vortex
cases will be different.

We are not sure what kind of monopole configurations lie in the various
spaces $\M_k^m$, but we conjecture that for $m = 0$ (or $m = k$), the
energy density is always confined to a finite neighbourhood of the plane
$x_3 = 0$, whereas for $0 < m < k$ it is possible for there to be monopole
clusters arbitrarily far from the plane $x_3 = 0$, arranged symmetrically
with respect to inversion in this plane. The examples discussed in Section
12 are consistent with this conjecture. If the roots of $q$ are distinct
and well-separated, then the configurations always consist of a set of
unit charge monopoles with their centres in the $x_3 = 0$ plane.

Our inversion formula is inconsistent with
Proposition 3.12 of \cite{AtiHit}.  There it was suggested that if
$$
R(z) = {p(z) \over q(z)} = \sum_i {\alpha_i \over z-\beta_i}
$$
is the rational map of a charge $k$ monopole, which consists of $k$
well-separated unit charge monopoles, then (using our conventions) the
individual monopoles are approximately located at the points
$(\beta_i, (1/2)\log|{\alpha_i}|)$
and have phases $\alpha_i|\alpha_i|^{-1}$.  Consideration
of our inversion formula suggests however that
the individual monopoles are
located at the points $(\beta_i, (1/2)\log|p(\beta_i)|)$ and have phases
$p(\beta_i) |p(\beta_i)|^{-1}$. This has recently
been proved by  Bielawski in \cite{Bie}. Interestingly the
spaces $\M_k^m$ play a distinguished role in Bielawski's work.

Finally notice that it follows from equations (\ref{eq:tau}) and
(\ref{eq:I}) and the fact that $\tau(\eta, 0) = I(\eta, 0)$ that
using $\tau(s)$ to construct the rational map is
the same as using $I(s)$, and hence the $p(\beta_i)$
occuring in the rational map defined using
$\tau(s)$ would be the reciprocal of the $p(\beta_i)$
we use, and would give the rational map as defined by
Hurtubise.

\section{Centred monopoles and rational maps}
\label{sec:centratmap}
We remarked earlier that although the positions and internal phases of the
$k$ \lq particles' in a charge $k$ monopole are only asymptotically
well-defined, every monopole has a well-defined centre and total phase.
This arises naturally in the twistor picture. If $S$ is the spectral curve
of a monopole then it intersects every fibre of $TP_1 \to P_1$ in $k$
points counted with multiplicity. If we add these points together we obtain
a new curve which is given by an equation $\eta + a_1(\zeta) = 0$. This
curve is a real section and hence $a_1$ is of the form
$$
a_1(\zeta) = -k((c_1 + ic_2) - 2c_3\zeta - (c_1 - ic_2)\zeta^2).
$$
The point $ c= (c_1, c_2, c_3) $ is the centre of the monopole.
To define the total phase requires a little more work.

Recall that the twistor space $Z$ for $S^1\times{\bf R}^3$ is the principal
bundle for the line bundle $L$ over $TP_1$. Its quotient $\bar Z$ by $\pm
1\in {\bf C}^*$ is the principal bundle for $L^2$. The spectral curve
$S\subset TP_1$ has a trivialization $s$ of $L^2$ and hence lifts to $\bar
Z$. Thus, over each point in $P_1$, the spectral curve defines $k$ points
counted with multiplicity in the fibre. This fibre is the group ${\bf
C}^*\times {\bf C}$ and we take the product of the points. This is a
section $s^k$ of $L^2$ over the curve $\eta+a_1(\zeta)=0$.

The bundle $L^2$ over any real section is trivial and
we fix as a choice of trivialisation $f$ over $\eta - k ((c_1 + i c_2) - 2
c_3 \zeta - (c_1 - i c_2)\zeta^2) = 0$
$$
\begin{array}{cc}
f_0(\eta, \zeta) =& \exp 2 k (c_3 + (c_1 - i c_2)\zeta)\\
f_1(\eta, \zeta) =& \exp 2 k (-c_3 + (c_1 + i c_2)/\zeta).
\end{array}
$$
It is easy to check that this non-vanishing section $f$
satisfies $\tau(f)f = 1$. Moreover
because $\tau(s)s = 1$ we must have $\tau(s^k) s^k = 1$.
 If we divide $s^k$ by $f$ we obtain
a holomorphic function which must be constant. In fact because $\tau(s^k)
s^k = 1$
and $\tau(f) f = 1$ this constant is a complex number of modulus $1$.
We define $s^k/f$ to be the total phase of the monopole. Notice that
if we act on the monopole by a constant diagonal gauge transformation
$\mu$ then $s$ is replaced by $\mu^2 s$ and the total phase is
multiplied by $\mu^{2k}$.

Let us now see how to construct the centre and total
phase of a monopole from its rational map. Notice first that if
we restrict the equation of the spectral curve to the fibre $\zeta = 0$
we  obtain an equation of the form
$$
\eta^k - k(c_1 + i c_2) \eta^{k-1} + ... = 0
$$
and hence $c_1 + ic_2$ is the average of the points of
intersection of the spectral curve with $\zeta = 0$ or
the average of the zeros of $q$.

Comparing the construction of the rational map
of a monopole we see that
$$
s^k_0(k(c_1+ic_2), 0) = \prod_i p(\beta_i) = \triangle(p,q)
$$
the resultant of $p$ and $q$.
It follows that
$$
{s^k\over f} = \triangle(p, q) \exp(-2kc_3).
$$
So:
\begin{prop} If $R(z) = p(z)/q(z)$ is the rational map of a monopole with
$q_0$ the average of the roots of $q$ and $\triangle(p,q)$ the resultant
of $p$ and $q$, then the centre of the monopole is
$$
(q_0, (1/2k)\log|\triangle(p,q)|)
$$
and the total phase is
$$
\triangle(p,q) |\triangle(p,q)|^{-1}.
$$
\end{prop}
It follows that  a monopole is centred if and only if the
zeroes of $q$ sum to zero and $|\triangle(p,q)| = 1$.
It will be useful to use a stronger notion of centring than this.
We call a monopole strongly centred if it is centred and the
total phase is $1$. From what we have just proven a monopole
is strongly centred if and only if its rational map satisfies
\begin{equation}
q_0 = 0 \quad\hbox{and}\quad \triangle(p, q) = 1.  \label{eq:strongcent}
\end{equation}

The resultant condition $\triangle(p, q) = 1$ was used in [\cite{AtiHit},
p.30] to identify the universal covering of the moduli space of centred
monopoles, but for a fixed complex structure in the hyperk\"ahler family.
Our description of strong centring gives an invariant approach, valid for
all complex structures. Considered in the context of the twistor space of
a hyperk\"ahler metric, it identifies the factor $X$ in the isometric
splitting [\cite{AtiHit}, p.34] $\tilde M_k = X \times  S^1\times\R^3 $ of
a $k$-fold covering of $M_k$ with the space of strongly centred monopoles.
It follows that the space of strongly centred monopoles is a geodesic
submanifold of $M_k$.
\vskip 0.2cm
\noindent{\bf Remark:}
Note that the twistor space $\bar Z$ is the twistor space for the
(trivial) hyperk\"ahler metric on the moduli space of 1-monopoles. A
$k$-monopole's centre and total phase then associates a 1-monopole
with a ${\bf Z}_k$-ambiguity of phases to a $k$-monopole.

\section{Symmetric curves in $TP_1$}
\label{sec:symmcur}
In eq.~(\ref{eq:spec}) we presented the general form of curves in $TP_1$
that occur as spectral curves of charge $k$ monopoles.  The coefficients
$a_r(\zeta)$ must satisfy the reality condition (\ref{eq:real}), and the
curve is centred at the origin in $\R ^3$ if $a_1(\zeta) = 0$. Here we
shall discuss the form of these curves when they are required to be
invariant under certain groups of rotations about the origin.

Let us recall that in $TP_1$, the $P_1$ of lines
through the origin are parametrized by $\zeta$ with $\eta = 0$.  The
line in the direction of the Cartesian unit
vector $(x_1, x_2, x_3)$ has
$\zeta = (x_1 + i x_2)/ (1+x_3)$.  It will be important to consider the
homogeneous coordinates $[\zeta_0, \zeta_1]$ on $P_1$, as well as the
inhomogeneous coordinate $\zeta = \zeta_1 /\zeta_0.$
An $SU(2)$ M\"obius transformation on the homogeneous coordinates,
$[\zeta_0, \zeta_1] \rightarrow [\zeta'_0, \zeta'_1]$, of the form
\begin{equation}
\begin{array}{rcl}
\zeta'_0 & = &- (b + ia) \zeta_1 + (d - ic)\zeta_0 \\
\zeta'_1 & = & (d + ic)\zeta_1 + (b - ia)\zeta_0   \\
\end{array} \label{eq:zetas}
\end{equation}
where $a^2 + b^2 + c^2 + d^2 =1$, corresponds to an $SO(3)$
rotation in ${\bf R} ^3$.  The rotation is by an angle $\theta$ about
the unit vector $(x_1, x_2, x_3)$, where
$
x_1 \sin {\theta \over 2} = a,\quad x_2 \sin {\theta \over 2} = b,
\quad x_3 \sin {\theta \over 2} = c, \quad
{\cos {\theta \over 2} = d}
$.
The inhomogeneous coordinate $\zeta$ transforms to
\begin{equation}
{\zeta' = {(d + ic)\zeta + (b - ia)
 \over -(b + ia)\zeta + (d-ic)}}.       \label{eq:zeta}
\end{equation}
Since $\eta$ is the coordinate in the tangent space to $P_1$ at $\zeta$,
it follows that if $\zeta$ transforms to $\zeta'$ as in (\ref{eq:zeta})
then $\eta$
transforms to $\eta'$ via the derivative of (\ref{eq:zeta}), that is
\begin{equation}
\eta' = {\eta\over (-(b + ia)\zeta + (d - ic))^2}. \label{eq:eta}
\end{equation}
A curve $P(\eta,\zeta)= 0$ in $TP_1$ is invariant under the M\"obius
transformation if $P(\eta',\zeta')=0$ is the same curve.  If the curve is
the spectral curve of a monopole, then the monopole is invariant under
the associated rotation.

The simplest group of symmetries is the cyclic group of rotations
about the $x_3$-axis, $C_n$.  The generator is the M\"obius
transformation
$$
\zeta' = e^{2\pi i/ n} \zeta, \quad \eta' = e^{2\pi i/ n}\eta.
$$
A curve $P(\eta,\zeta) = 0$ is invariant if all terms of $P$ have the same
degree, mod $n$.  A curve of the form (\ref{eq:spec}) is $C_n$-invariant
if all terms have degree $k$, mod $n$. In particular, it is
$C_k$-invariant if all terms have degree zero, mod $k$.

For there to be axial symmetry about the $x_3$-axis, with symmetry group
$C_\infty$, the curve must be invariant under $\zeta \rightarrow
e^{i\theta} \zeta, \quad \eta\rightarrow e^{i \theta} \eta$, for all
$\theta$. This requires that all terms in $P(\eta,\zeta)$ have degree $k$.
There is a unique axially symmetric, strongly centred monopole for each
charge $k$. It is shown in \cite{Hit} that its spectral curve is
$$
\begin{array}{cl}
\eta  \prod^{m}_{l=1} (\eta^2 + l^2 \pi^2 \zeta^2) & = 0 \quad \hbox
{for}\quad k=2m +1
\\
\prod^{m}_{l=0} \left (\eta^2 + (l+ {1 \over 2})^2\pi^2 \zeta^2 \right )
& = 0 \quad \hbox {for}\quad k = 2m + 2.
\end{array}
$$
Notice that these curves are not determined by symmetry alone, and that the
coefficients of $P$ are transcendental numbers.  The only curve of the
form (\ref{eq:spec}) which has full $SO(3)$ symmetry is $\eta^k = 0$.
This is the spectral curve of a unit charge monopole at the origin when
$k=1$, but for $k > 1$ it is not the spectral curve of a monopole.

The groups $C_n$ and $C_\infty$ are extended to the dihedral groups $D_n$
and $D_\infty$ by adding a rotation by $\pi$ about the $x_1$-axis.  This
rotation corresponds to the transformation on $TP_1$
$$
\zeta' = {1 \over \zeta}, \quad \eta' = - {\eta \over \zeta^2}.
$$
Under this transformation, and for any constant $\nu$,
$$
(\eta^2 + \nu \zeta^2)' = {1 \over \zeta^4} (\eta^2 + \nu \zeta^2),
$$
so each of the axially symmetric monopoles has symmetry group
$D_\infty$.

Recall from Section \ref{sec:invratmap} that $P(\eta, \zeta) = 0$ is
reflection symmetric under $x_3 \rightarrow -x_3$ if $P$ is even in
$\zeta$. By a similar argument to that in Section \ref{sec:invratmap}, the
reflection $x_2 \rightarrow -x_2$ corresponds to $\zeta \rightarrow
\overline\zeta , \eta \rightarrow \overline\eta$, so a curve $P(\eta,
\zeta) = 0$ is invariant under this reflection if all coefficients in
$P(\eta, \zeta)$ are real. The axially symmetric monopoles therefore have
these reflection symmetries too.

As an example of finite cyclic or dihedral symmetry, let us consider centred
$k=3$ curves with either $C_3$ or $D_3$ symmetry.  Before imposing the
symmetry, the curves are of the form
\begin{eqnarray}
\eta^3  &+ &\eta(\alpha_4 \zeta^4 + \alpha_3 \zeta^3 + \alpha_2 \zeta^2 +
\alpha_1 \zeta + \alpha_0) \nonumber \\
   & & \mbox{}+ (\beta_6 \zeta^6+ \beta_5 \zeta^5 + \beta_4 \zeta^4 +
\beta_3 \zeta^3 +
\beta_2 \zeta^2 + \beta_1 \zeta + \beta_0) = 0   \label{eq:symmk=3}
\end{eqnarray}
subject to the reality conditions
$$
\begin{array}{cl}
\alpha_4 & = \overline{\alpha}_0, \ \alpha_3  =-\overline{\alpha}_1, \
\alpha_2 =
\overline{\alpha}_2, \\
\beta_6 &  = -\overline{\beta}_0,\ \beta_{5} = \overline{\beta}_1, \
\beta_{4} =
-\overline{\beta}_2, \ \beta_3 = \overline{\beta}_3.

\end{array}
$$
$C_3$ symmetry implies that (\ref{eq:symmk=3}) reduces to
\begin{equation}
\eta^3 + \alpha \eta \zeta^2 + \beta \zeta^6 +
\gamma \zeta^3 - \bar{\beta} = 0  \label{eq:C3symm}
\end{equation}
where $\alpha$ and $\gamma$ are real.  By a rotation about
the $x_3$-axis, we can orient the curve so that $\beta$ is real, too,
and then there is reflection symmetry under $x_2 \rightarrow -x_2$.
There is $D_3$ symmetry if $\gamma = 0$; then the curve reduces to
\begin{equation}
\eta^3 + \alpha \eta \zeta^2 + \beta(\zeta^6 -1) = 0  \label{eq:D3symm}
\end{equation}
with $\alpha$ and $\beta$ real.

The axisymmetric charge $3$ monopole has a spectral curve of type
(\ref{eq:D3symm})
with $\alpha = \pi^2$ and $\beta = 0$.  Also, three
well-separated unit charge monopoles at the vertices of an equilateral
triangle can have $D_3$ symmetry. The
spectral curve is asymptotic to the product of three stars at
$$
(x_1, x_2, x_3) =\left \{ (a, 0, 0), \quad (a \cos
{2\pi\over 3}, a \sin {2\pi\over 3}, 0),
\quad (a \cos {4\pi \over 3}, a \sin
{4\pi \over 3}, 0) \right
\},
$$
that is,
$$
(\eta - a (1 - \zeta^2))(\eta - a\omega (1 - \omega \zeta^2)) (\eta
- a\omega^2 (1 - \omega^2 \zeta^2) ) = 0,
$$
where $\omega = e^{2 \pi i/3}.$  Multiplied out, this is a curve of the
form (\ref{eq:D3symm}) with $\alpha = 3a^2$ and $\beta =a^3$, or
equivalently $\alpha^3 = 27\beta^2$. We shall find out more about the
spectral curves of charge 3 monopoles with symmetry $C_3$ or $D_3$ when we
consider the rational maps associated with the monopoles (see Section 12).

Symmetry under $C_4$ is rather a weak constraint on curves with $k = 4$.
On the other hand $D_4$ symmetry implies that a $k = 4$ curve is of the
form \begin{equation} \eta^4 + \alpha \eta^2 \zeta^2 + \beta \zeta^8 +
\gamma\zeta^4 + \beta = 0 \label{eq:D4symm} \end{equation} with $\alpha,\
\beta$ and $\gamma$ real. The axisymmetric charge $4$ monopole has this
form of spectral curve, with $\alpha = (5/2) \pi^2$, $\beta = 0$ and
$\gamma = (9/16) \pi^4$.  Four well-separated unit charge monopoles at the
vertices of the square $\left \{ (\pm a, 0, 0), (0, \pm a, 0) \right \}$
can have $D_4$ symmetry. The spectral curve is asymptotic to a product of
stars, and is of the form (\ref{eq:D4symm}), with $\alpha = 4a^2,\ \beta =
- a^4$ and $\gamma = 2a^4$. After a $\pi/4$ rotation, the monopoles are at
$(\pm a/\sqrt{2}, \pm a/\sqrt{2}, 0)$, and $\alpha = 4a^2, \beta = a^4$
and $\gamma =2a^4$.

There is another interesting asymptotic 4-monopole configuration, with a
spectral curve of type (\ref{eq:D4symm}).  Consider two well-separated
axisymmetric
charge $2$ monopoles, centred at $(0, 0, b)$ and $(0, 0, -b)$, and with
the $x_3$-axis the axis of symmetry.  The spectral curve is asymptotic
to a product of curves associated with the charge $2$ monopoles.
The spectral curve of a centred axisymmetric charge $2$
monopole is $\eta^2 + {1 \over 4} \pi^2 \zeta^2 = 0$.  This factorizes
as $(\eta + {1 \over 2} i\pi\zeta) (\eta - {1 \over 2} i\pi \zeta) = 0$,
which is a product of stars at the complex
conjugate points $(0,0, \pm i\pi/4).$  Translation by $b$ gives
the curve
$$
\eta^2 + 4b \eta \zeta + (4b^2 + {1 \over 4} \pi^2) \zeta^2 = 0
$$
which is the product of stars at $(0,0, b \pm i\pi/4)$.  Similarly,
translation by $-b$ gives
$$
\eta^2 - 4b \eta \zeta + (4b^2 + {1 \over 4} \pi^2) \zeta^2 =0
$$
and the product of these is the curve
$$
\eta^4 + ({1 \over 2} \pi^2 - 8b^2) \eta^2 \zeta^2 +
(4b^2 + {1 \over 4} \pi^2)^2 \zeta^4 = 0.
$$
Since all terms have degree 4 this curve is axisymmetric; however, the
actual spectral curve of the charge $4$ monopole has symmetry $D_4$, as we
shall see in Section \ref{sec:scatt}, becoming axisymmetric only in the
limit of infinite separation.

Let us now investigate the curves in $TP_1$ with the symmetries of a
regular solid. Some of these are special cases of the curves we have
already discussed.  There  are three rotational symmetry groups to
consider, those of a tetrahedron, an octahedron and an icosahedron. The
direct way to construct a symmetric curve is to find M\"obius
transformations which generate the symmetry group, and calculate the
conditions for the curve to be invariant under all of them.  For example,
a curve of type (\ref{eq:D4symm}), with $D_4$ symmetry, has octahedral
symmetry if it is invariant under the transformation
$$
\zeta' = {i \zeta + 1 \over \zeta + i} , \quad \eta' =
{-2 \over (\zeta + i)^2 }\eta,
$$
which corresponds to a $\pi/2$ rotation about the $x_1$-axis, and this
requires that the curve reduces to
$$
\eta^4 + \beta (\zeta^8 + 14 \zeta^4 + 1) = 0.
$$
A more powerful and less laborious approach is to use the theory of
invariant bilinear forms and polynomials on $P_1$, as expounded in Klein's
famous book \cite{Kle}.

Consider a homogeneous bilinear form $Q_r(\zeta_0, \zeta_1)$ of degree
$r$, and its associated inhomogeneous polynomial $q_r(\zeta)$ defined by
$$
Q_r(\zeta_0, \zeta_1) = \zeta^r_0 q_r (\zeta).
$$
Generally $q_r$ has degree $r$, but it may have lower degree.  Suppose
$Q_r (\zeta_0, \zeta_1)$ is invariant under a M\"obius transformation of
the form (\ref{eq:zetas}).  Then $q_r(\zeta)$ transforms in a simple
way under the corresponding transformation (\ref{eq:zeta}), namely
\begin{equation}
q'_r (\zeta) = {q_r(\zeta) \over
\left (-(b + ia) \zeta + (d - ic)\right )^r}.  \label{eq:qtrans}
\end{equation}
On the other hand, $\eta$ transforms as in (\ref{eq:eta}).  Consider a
centred curve in $TP_1$,
$$
P(\eta, \zeta) \equiv \eta^k + \eta^{k-2} q_4
(\zeta) + \eta^{k-3} q_6 (\zeta) +
\dots +
q_{2k} (\zeta) = 0.
$$
If, under a M\"obius transformation, each polynomial $q_r(\zeta)$
transforms as in (\ref{eq:qtrans}), and $\eta$ as in (\ref{eq:eta}), then
each term in the polynomial $P(\eta, \zeta)$ is multiplied by the same
factor $\left (-(b + ia) \zeta + (d - ic) \right )^{-2k}$, so the curve is
invariant.  It follows that curves invariant under the rotational symmetry
group of a regular solid can be constructed from the inhomogeneous
polynomials $q_r$ derived from the bilinear forms $Q_r$ invariant under
the group.

Let $G$ denote the tetrahedral, octahedral or icosahedral group. Klein has
described the ring of bilinear forms, $Inv_G$, which change only by a
constant factor under each transformation of $G$ $\, - \,$ for each form
these factors define an abelian character of $G$.  Let $Inv_G^\star$ be
the subring of strictly invariant forms.  A form $Q$ is in $Inv_G$ if the
roots of the associated polynomial $q$ are invariant under $G$, that is,
if they are the union of a set of $G$-orbits on $P_1$.

Generic $G$-orbits on $P_1$ consist of $|G|$ points, i.e. 12, 24 and
60 points respectively for the three groups.  The associated forms of
degree $|G|$ are always strictly invariant under $G$, and they span a
vector space of forms, of dimension two.  For each group $G$, there
are also three forms of degree less than $|G|$ associated with special
orbits of $G$, and these generate the ring $Inv_G$.  Let $V$, $E$ and
$F$ be the set of vertices, mid-points of edges, and centres of faces
of the centred regular solid (tetrahedron, octahedron or icosahedron)
invariant under $G$.  Centrally project these points onto the unit
sphere, identified with $P_1$, denoting the resulting sets of points
again by $V, E$ and $F$.  $V$ is a $G$-orbit, so there is a form $Q_V$
in $Inv_G$ and an associated polynomial $q_V$, such that $Q_V$ has
degree $|V|$ and $Q_V = 0$ at all points of $V$.  Similarly, there are
forms and polynomials $Q_E, Q_F$ and $q_E,q_F$.  Table 1 gives the
polynomials $q_V, q_E$ and $q_F$ for the three groups $G$, and a star
indicates that the associated form $(Q_V, Q_E$ or $Q_F)$ is strictly
$G$-invariant.  [A choice of orientation has been made for the solids:
the tetrahedron has its vertices at $(1/\sqrt{3}) (\pm 1, \pm 1, \pm
1)$, with either two or no signs negative; the octahedron has its vertices
on the Cartesian axes; the icosahedron has two vertices on the
$x_3$-axis and is invariant under the dihedral group $D_5$.]

All the icosahedral forms are strictly invariant because the
icosahedral group $A_5$ is simple, and has no non-trivial abelian
characters.  The tetrahedral forms $Q_V$ and $Q_F$ are not strictly
invariant, but acquire factors of $e^{\pm 2\pi i/3}$ under a $2\pi/3$
rotation about a 3-fold symmetry axis; so $Q_V Q_F$ is strictly invariant.
In fact, the polynomial associated with $Q_V Q_F$ is $\zeta^8 + 14\zeta^4
+ 1$, which has octahedral symmetry.  Similarly, the octahedral forms
$Q_V$ and $Q_E$ acquire factors of $-1$ under a rotation by $\pi/2$
around a 4-fold symmetry axis, and $Q_V Q_E$ is strictly invariant.

There are remarkable identities satisfied by the forms  $Q_{V}, Q_{E}$ and
$Q_F$ (which remain true if the forms $Q$ are replaced by the
associated polynomials $q$), namely
$$
\begin{array}{cl}
Q_V^3 - Q^3_F - 12\sqrt{3}i \ Q^2_E & = 0 \quad \hbox
{for the tetrahedral group}
\cr
108 \ Q_V^4 -  Q^3_F + Q_E^2 & = 0  \quad \hbox {for the octahedral group} \\
1728 \ Q^5_V - Q_F^3 - Q^2 _E & = 0 \quad \hbox
{for the icosahedral group}.
\end{array}
$$
These identities occur, because each term is a strictly invariant form
of degree $|G|$, lying in the two-dimensional vector space of forms
associated with the generic $G$-orbits.

% THIS IS THE BEGINNING OF TABLE 1
$$\vbox{
\offinterlineskip
\halign{
\strut\vrule \hfil $\;$#$\;$\hfil & \vrule \hfil $\;#\;$
\hfil & \vrule \hfil
$\;#\;$ \hfil & \vrule \hfil $\;#\;$ \hfil \vrule \cr
\noalign{\hrule}
&&& \cr
G & q_{_V} & q_{_E} & q_{_F} \cr
&&& \cr
\noalign{\hrule}
&&& \cr
Tetrahedral & \zeta^4 + 2\sqrt{3} i \zeta^2 + 1& \zeta(\zeta^4 -1)^\star &
\zeta^4 -2\sqrt{3} i \zeta^2 + 1 \cr
&&& \cr
\noalign{\hrule}
&&& \cr
Octahedral & \zeta(\zeta^4 -1) & \zeta^{12} - 33\zeta^8 & \zeta^8 +
14\zeta^4 + 1^\star\cr
&& -33\zeta^4 + 1 & \cr
&&& \cr
\noalign{\hrule}
&&& \cr
Icosahedral & \zeta(\zeta^{10} + 11 \zeta^5 -1)^\star & \zeta^{30} + 522
\zeta^{25}
& \zeta^{20} - 228 \zeta^{15} + 494 \zeta^{10} \cr
& & -10005 \zeta^{20} - 10005 \zeta^{10} &  \cr
&& -522 \zeta^5+ 1^\star & + 228\zeta^5 + 1^\star  \cr
&&& \cr
\noalign{\hrule}
}
}$$

\vskip 5pt
\noindent
{\narrower\smallskip\noindent Polynomials associated with the special
orbits $V, E$ and $F$ of the rotational symmetry groups of the regular
solids.  A star $(\star)$ denotes that the homogeneous bilinear form
$Q$ related to the polynomial $q$ is strictly invariant.\smallskip}
\medskip
\begin{center}
{\bf Table 1}
\end{center}
\bigskip

% THIS IS THE END OF TABLE 1

We can now write down some examples of invariant curves in $TP_1$,
also satisfying the reality conditions (\ref{eq:real}).
Recall that invariant curves in $TP_1$ must be constructed from
polynomials derived from strictly invariant forms.  The simplest curves
with tetrahedral symmetry are
$$
\eta^3 + i a \zeta(\zeta^4 - 1) = 0
$$
where $a$ is real.  After a rotation, these become
$$
\eta^3 + a (\zeta^6 + 5\sqrt{2}\zeta^3 - 1) = 0,
$$
which are of the form (\ref{eq:C3symm}), exhibiting manifest $C_3$ symmetry
about the $x_3$-axis. The simplest curves in $TP_1$ with octahedral
symmetry are
$$
\eta^4 + a(\zeta^8 + 14 \zeta^4 + 1) = 0  \label{eq:octa}
$$
with $a$ real. More generally, the $k=4$ curves with tetrahedral
symmetry are of the form
$$
\eta^4 + ib\eta\zeta(\zeta^4 - 1) + a(\zeta^8 + 14 \zeta^4 + 1) = 0
\label{4tetra}
$$
with $a$ and $b$ real.  Finally, the simplest curves
with icosahedral symmetry are
$$
\eta^6 + a\zeta(\zeta^{10} + 11 \zeta^5 - 1) = 0 \label{eq:icosa}
$$
with $a$ real.

We shall discuss in the next Sections the possibility that some of
these curves are spectral curves of monopoles.

\section{Symmetric curves and elliptic curves}
\label{sec:symmellip}
Let $G\subset SO(3)$ be the symmetry group of a regular solid and $\tilde
G\subset SU(2)$ the corresponding binary group. Let $V$ be the
2-dimensional defining representation of $SU(2)$. The $n$-th symmetric
power $S^nV$ may be considered as the action of $SU(2)$ on homogeneous
polynomials of degree $n$ in $\zeta_0,\zeta_1$. Alternatively, the
representation is on the space of holomorphic sections of the line bundle
${\cal  O}(n)$  on $P_1$, which is described as the polynomials of degree
$\le n$ in the affine parameter $\zeta=\zeta_1/\zeta_0$.

Let $n$ be the smallest degree of a homogeneous polynomial invariant under
$\tilde G$. For the  tetrahedral, octahedral and icosahedral group the
value of $n$ is $2k$ where  $k$ is  respectively  3,4,6 and in each case
there is a unique (up to a multiple) invariant  polynomial as given in
Table 1. Note that in all cases $\vert G \vert =2k(k-1)$. In the previous
Section, we defined curves in $TP_1$ by
 \begin{eqnarray}
 \eta^3+ia\zeta(\zeta^4-1)&=&0 \label{A}\\
 \eta^4+a(\zeta^8+14\zeta^4+1)&=&0 \label{B}\\
 \eta^6+a\zeta(\zeta^{10}+11\zeta^5-1)&=&0\label{C}
 \end{eqnarray}
where $a$ is real. Each such curve $S$ is invariant by the appropriate
group $G$ and satisfies the reality conditions for a spectral curve. We
shall prove that the first two are indeed spectral curves for non-singular
monopoles of charges $3$ and $4$, for suitable values of $a$. On the other
hand, by the same methods we shall see that there is no charge 6 monopole
which has icosahedral symmetry.

The method consists of finding explicitly the solution to Nahm's equations
with these spectral curves,  thereby solving the non-linear part of the
monopole problem. To find the monopole configuration itself involves
solving an associated linear differential equation \cite{AtiHit}.

The solution to Nahm's equations will be expressed in terms of elliptic
functions, and the constant $a$ in terms of the periods of an elliptic
curve. The advance evidence for this fact lies in the following result:
\begin{prop}
The curve $S$ in (\ref{A}),(\ref{B}),(\ref{C}) is  smooth and of genus
$(k-1)^2$. Its quotient by $G$ is an elliptic curve.
\end{prop}
If we think of the solution to Nahm's equations given in \cite{Hit1} as a
linearization of the equations on the Jacobian of the spectral curve, then
the above proposition implies that a $G-$invariant solution is linearized
on the Jacobian of the elliptic curve and hence expressible in terms of
elliptic functions. We shall achieve this directly, however, without
making use of the general method of \cite{Hit1}. \vskip .5cm \noindent{\it
Proof of Proposition}: Smoothness is immediate since the polynomials by
inspection have distinct roots. The genus formula is standard
\cite{Hit1,AtiHit}.

Now consider the action of $G$ on $S$. Let $m$ be a  fixed point of an
element of $G$. Consider its image in $P_1$. The group $G$ acts on the
spectral curve through the natural action on the tangent bundle $TP_1$,
but for isometries the action of an isotropy group on the tangent space of
the point is faithful. Thus the zero vector is the only one fixed. But
this is given by $\eta=0$ and from the equation of the curve these points
are the $2k$ zeros of the polynomial.

It follows that  $G$ acts freely except at the $2k$ points $\eta=0$, and
since $\vert G \vert=2k(k-1)$, the stabilizer of each is  of order
$(k-1)$. These are the stabilizers of vertices of the regular solids, and
hence cyclic. Thus by the Riemann-Hurwitz formula, the genus $g$ of the
quotient satisfies $$2-2(k-1)^2=2k(k-1)(2-2g)-2k(k-2)$$ and hence $g=1$.

\vskip 1cm
\section{Special solutions to Nahm's equations}
\label{sec:Nahm}

Recall that a centred charge $k$ monopole may be obtained from a
solution of Nahm's equations \cite{Hit1}
\begin{eqnarray}
\frac{dT_1}{ds}=[T_2,T_3]\nonumber\\
\frac{dT_2}{ds}=[T_3,T_1]\label{N}\\
\frac{dT_3}{ds}=[T_1,T_2]\nonumber
\end{eqnarray}
where $T_i(s)$ is a function with values in the Lie algebra ${\frak
{su}}(k)$. It must satisfy moreover a reality condition
$T_i(2-s)=T_i(s)^T$ with respect to an orthonormal basis compatible with
the unitary structure. (In the explicit formulae which follow, we have not
always used such a basis, preferring one which is simpler for
calculations. We rely on the reality of the spectral curve and the
description of \cite{Hit1} to assure the existence of such a  basis.) The
solution to Nahm's equations must be regular for $s\in(0,2)$ and have
simple poles at $s=0$ and $s=2$. At $s=0$, we have an expansion
$$T_i=R_i/s+\dots .$$
Nahm's equations themselves imply that, at {\it any} simple pole, the
residues satisfy
\begin{eqnarray*}
R_1=-[R_2,R_3]\\
R_2=-[R_3,R_1]\\
R_3=-[R_1,R_2]
\end{eqnarray*}
and so define a representation of the Lie algebra ${\frak{so}}(3)$. In
order for a solution to give
a monopole, the representation at $s=0$ and $s=2$ must
(see \cite{Hit1} or \cite{AtiHit} Chapter 16) be the unique irreducible
representation $S^{k-1}V$ of dimension $k$.
In fact, as shown
in \cite{Hit1}, this space is canonically isomorphic to $H^0(P_1,{\cal
O}(k-1))$ under the projection from the spectral curve $S$.
With any solution to Nahm's equations, the coefficients of the polynomial
\begin{equation}
P(\eta,\zeta) = \det(\eta +i(T_1+iT_2)-2iT_3\zeta -i(T_1-iT_2)\zeta^2)
\label{eq:Nahmpoly}
\end{equation}
are independent of $s$, and indeed for a monopole $P(\eta,\zeta) = 0$
defines $S$.

We may regard the triple of Nahm matrices as
a function
with values in
$${\bf R}^3\otimes {\frak {su}}(k).$$
The action of the rotation group $SO(3)$ on a monopole then appears as the
tensor product of its natural action on ${\bf R}^3$ and ${\frak {su}}(k)$.
In terms of the irreducibles $S^mV$, this is the representation
$$S^2V\otimes \End_0 (S^{k-1}V)$$ where $\End_0$ denotes trace zero
endomorphisms.

In the above situation where the monopole is $G$-invariant, the Nahm
matrices lie in a  subspace of $S^2V\otimes \End_0(S^{k-1}V)$ which is
fixed by $G$.
\vskip .5cm
\begin{prop} The fixed point set of $G$ in $S^2V\otimes \End_0(S^{k-1}V)$
is a 2-dimensional vector space.
\end{prop}
\noindent{\it Proof:} First take the Clebsch-Gordan decomposition of
$\End(S^{k-1}V)\cong S^{k-1}V\otimes S^{k-1}V$ into irreducibles of
$SO(3)$. We obtain
$$
\End_0(S^{k-1}V)\cong S^{2k-2}V\oplus S^{2k-4}V\oplus \dots \oplus S^2V.
$$
The single $S^2V$ factor gives a 1-dimensional subspace of $S^2V\otimes
S^{k-1}V\otimes S^{k-1}V$ fixed by $SO(3)$. This is simply the
homomorphism of Lie algebras $${\frak {so}}(3)\cong S^2V \rightarrow
\End_0(S^{k-1}V)$$ defined by the representation $S^{k-1}V$. Now the
Clebsch-Gordan decomposition also gives $$S^{2m}V\otimes
S^2V=S^{2m+2}V\oplus S^{2m}V \oplus S^{2m-2}V.$$ But  $2k$ is, by choice,
the smallest positive integer $n$ such that $S^nV$ has a $G$-invariant
vector. Thus $S^{2m}V\otimes S^2V$ has no invariants for $1<m<k-1$, and
for $m=k-1$ there is a unique one lying in $S^{2k}V$. This gives another
1-dimensional fixed subspace, and therefore a 2-dimensional space
altogether.
\vskip 0.5cm

We use this fact next to give a substantial simplification of Nahm's
equations in the $G$-invariant case. From Proposition 5, any $G$-invariant
element $T=(T_1,T_2,T_3)$ of ${\bf R}^3\otimes{\frak su}(k)$ can be
written as
$$T_i=x\rho_i+yS_i$$
where $\rho:{\bf R}^3 \rightarrow {\frak {su}}(k)$ is the representation of
${\frak{so}}(3)$ on ${\bf C}^k$ and $(S_1,S_2,S_3)$ is the $G$-invariant
vector in $S^{2k}V\subset {\bf R}^3\otimes {\frak {su}}(k)$. In particular
the Nahm matrices for a $G$-invariant monopole can be expressed as
\begin{equation}
T_i(s)=x(s)\rho_i +y(s)S_i \qquad {\rm for} \quad i=1,2,3.
\label{Nahm}
\end{equation}
Given one invariant element $T$ of ${\bf R}^3\otimes {\frak {su}}(k)$, we
can use the cross product on ${\bf R}^3$ and the Lie bracket on
$\End_0({\bf C}^k)$ to define another, $T\times T$. Since this again lies
in the two-dimensional fixed subspace generated by $\rho$ and $S$, there
must be constants $\alpha,\beta,\gamma,\delta$ such that:
$$[S_1,\rho_2]+[\rho_1,S_2]=\alpha \rho_3 + \beta S_3$$
$$[S_1,S_2]= \gamma \rho_3 + \delta S_3$$
and corresponding expressions obtained by cyclic permutation. The analogous
term for $\rho_i$ is determined by the fact that it gives a representation
of ${\frak{so}}(3)$:
$$[\rho_1,\rho_2]=2\rho_3 \qquad {\rm etc.}$$
using the standard basis of the Lie algebra. From this, Nahm's equations
(\ref{N}) become
\begin{eqnarray}
\frac{dx}{ds}&=&2x^2+\alpha xy+\gamma y^2 \label{e1}\\
\frac{dy}{ds}&=&\beta xy + \delta y^2 . \label{e2}
\end{eqnarray}
Equations in this general form can always be reduced to quadratures, but in
our case we shall calculate the precise values of the constants $\alpha,
\beta, \gamma$ and $\delta$ and find $x(s)$ and $y(s)$ exactly.

To evaluate $\alpha,\beta,\gamma$ and $\delta$ we must find the $k\times k$
matrices $\rho_i$ and $S_i$. If $e_1,e_2,e_3$ is the standard basis of the
Lie algebra ${\frak{so}}(3)$, with $[e_1,e_2]=2e_3$ etc. then
$\rho_i=\rho(e_i)$ where
$$\rho:{\frak{so}}(3)\rightarrow \End_0({\bf C}^k)$$
is the irreducible $k$-dimensional representation. If we regard this as the
action on homogeneous polynomials
$$a_0\zeta_0^{k-1}+a_1\zeta_0^{k-2}\zeta_1+\dots+a_{k-2}\zeta_0\zeta_1^{k-2}
+a_{k-1}\zeta_1^{k-1}$$
then setting
$$X=\frac{1}{2}(e_1-ie_2),\quad Y=-\frac{1}{2}(e_1+ie_2),\quad H=-ie_3$$
we have the Lie brackets
\begin{equation}
[X,Y]=H,\quad [H,X]=2X, \quad [H,Y]=-2Y
\label{XYH}
\end{equation}
and the representation is defined on polynomials by the operators
$$X=\zeta_1\frac{\partial}{\partial \zeta_0},\quad Y=\zeta_0
\frac{\partial}{\partial \zeta_1},\quad H=-\zeta_0\frac{\partial}{\partial
\zeta_0}+ \zeta_1\frac{\partial}{\partial \zeta_1}.$$

To determine the $S_i$, we have to reinterpret each of the
polynomials in $\zeta$ in (\ref{A}),(\ref{B}),(\ref{C}), using the inclusions
$$
S^{2k}V\subset S^{2}V\otimes S^{2k-2}V\subset S^2V\otimes \End_0(
S^{k-1}V).
$$
The first inclusion is simply polarization (or differentiation) of a
homogenous polynomial $Q_{2k}(\zeta_0,\zeta_1)$ of degree $2k$
$$
Q_{2k}\mapsto \xi_0^2\otimes \frac{\partial^2Q_{2k}}{\partial
\zeta_0^2}+2\xi_0\xi_1\otimes
\frac{\partial^2Q_{2k}}{\partial\zeta_0\partial\zeta_1}+\xi_1 ^2\otimes
\frac{\partial^2Q_{2k}}{\partial \zeta_1^2}.
$$
The second inclusion comes from
$S^{2k-2}V\subset \End_0(S^{k-1}V)$. A useful
way to view this is to see the image $\rho({\frak{so}}(3))$ as the principal
3-dimensional subalgebra of ${\frak {sl}}(k,{\bf C})$ (see \cite{Kos}).

Any complex simple Lie algebra ${\frak g}$ of rank $r$ breaks up into $r$
irreducible representations under the action of its principal
3-dimensional subalgebra $<X,Y,H>$. Moreover, the nilpotent element $X$ is
a regular nilpotent in ${\frak g}$, and belongs to an $r$-dimensional
abelian nilpotent subalgebra. The weight spaces of this algebra under the
action of $H$ are the highest weight spaces of the irreducible
representations into which ${\frak g}$ breaks up.

This is all for  a general Lie algebra. In our case, for ${\frak
{sl}}(k,{\bf C})$, the decomposition is into representation spaces
 $$
 S^{2k-2}V\oplus S^{2k-4}V\oplus \dots\oplus S^2V
 $$
so that the subspace $S^{2k-2}V$ is the representation which contains the
vector of highest weight among all elements in the Lie algebra. The
nilpotent element $X$ lies in the 3-dimensional Lie algebra $S^2V$, and,
being regular, has rank $k-1$. It acts cyclically on ${\bf C}^k$ and so
its centralizer is spanned by the powers of $X$. The  element  of  highest
 weight  which  commutes  with  it  is thus $X^{k-1}$,  of rank 1. This is
the highest weight vector of $S^{2k-2}V$ and so applying  $Y$  to  this
element we generate  the whole subspace $S^{2k-2}V\subset
\End_0(S^{k-1}V)$. Thus, given a homogeneous polynomial
$Q(\zeta_0,\zeta_1)\in S^{2k-2}V$, we define a $k\times k$ matrix $S(Q)$
by finding the polynomial $\tilde Q$ such that
$$
Q(\zeta_0,\zeta_1)={\tilde Q}(\zeta_0\frac{\partial}{\partial\zeta_1})
\zeta_1^{2k-2}
$$
and then setting
$$
S(Q)={\tilde Q}(\ad Y)X^{k-1}.
$$
In this way we evaluate the matrices $S_1,S_2,S_3$ for the three cases in
 (\ref{A}),(\ref{B}),(\ref{C}).

\section{Tetrahedral symmetry}
\label{sec:tetra}
In the tetrahedral case, $k=3$ and the irreducible representation $\rho$ on
${\bf C}^3$ is the adjoint representation. From (\ref{XYH}), relative to
the basis $X,Y,H$, the action of $X$ is given by the matrix
$$\left(
\matrix {0&-2&0\cr
            0&0&1\cr
            0&0&0\cr}
            \right)
            $$
 and then a highest weight vector is a multiple of
 $$X^2=\left(
\matrix {0&0&-2\cr
            0&0&0\cr
            0&0&0\cr}
            \right).$$
 The matrix $Y$ in this representation is
 $$Y= \left(
\matrix {0&0&0\cr
            -1&0&0\cr
            0&2&0\cr}
            \right)
            $$
and so, using the above procedure, we can evaluate $\rho_i$ and $S_i$. For
the representation $\rho$ we use a basis such that
 $$\rho_1=\left(
 \matrix{0&2&0\cr
         -1&0&2\cr
         0&-1&0\cr}
         \right),  \quad \rho_2=\left(
         \matrix{0&2i&0\cr
                 i&0&2i\cr
                 0&i&0\cr}
                 \right), \quad \rho_3=\left(
                 \matrix{2i&0&0\cr
                         0&0&0\cr
                         0&0&-2i}
                         \right).
                         $$
To find $S_i$ we first take the polynomial
$$
Q_6(\zeta_0,\zeta_1)=\zeta_0\zeta_1(\zeta_1^4-\zeta_0^4)
$$
and polarize it to obtain
$$
-20\xi_0^2\otimes
\zeta_0^3\zeta_1+10\xi_0\xi_1\otimes
(\zeta_1^4-\zeta_0^4)+20\xi_1^2\otimes \zeta_0\zeta_1^3.
$$
Relating the basis $\xi_0^2,\xi_0\xi_1,\xi_1^2$ to the standard orthonormal
basis of $S^2V\cong{\frak{so}}(3)$, we have $(S_1,S_2,S_3)$ given as a
multiple of
$$
(2i(\zeta_0\zeta_1^3-\zeta_0^3\zeta_1),
2(\zeta_0\zeta_1^3+\zeta_0^3\zeta_1
),(\zeta_1^4-\zeta_0^4))
$$
where we now have to interpret the quartic polynomials as $3\times
3$ matrices. Following the procedure described above, we obtain
$$S_1=\left(
\matrix{0&\frac{1}{2}&0\cr
        \frac{1}{4}&0&-\frac{1}{2}\cr
        0&-\frac{1}{4}&0\cr}
        \right), \quad S_2=i\left(
        \matrix{0&-\frac{1}{2}&0\cr
               \frac{1}{4}&0&\frac{1}{2}\cr
               0&-\frac{1}{4}&0\cr}
               \right), \quad S_3=i\left(
               \matrix{0&0&-1\cr
                      0&0&0\cr
                      \frac{1}{4}&0&0\cr}
                      \right).$$

\noindent{\bf Remark}: The above form is a direct consequence of
choosing a convenient basis for calculations, but
the   reality  conditions are more evident if we change basis.
We then find
$$\rho_1=\left(
 \matrix{0&0&0\cr
         0&0&-2\cr
         0&2&0\cr}
         \right),  \quad \rho_2=\left(
         \matrix{0&0&2\cr
                 0&0&0\cr
                 -2&0&0\cr}
                 \right), \quad \rho_3=\left(
                 \matrix{0&-2&0\cr
                         2&0&0\cr
                         0&0&0}
                         \right);
$$
 $$
 S_1=\left(
\matrix{0&0&0\cr
        0&0&-\frac{1}{2}\cr
        0&-\frac{1}{2}&0\cr}
        \right), \quad S_2=\left(
        \matrix{0&0&-\frac{1}{2}\cr
               0&0&0\cr
               -\frac{1}{2}&0&0\cr}
               \right), \quad S_3=\left(
               \matrix{0&-\frac{1}{2}&0\cr
                      -\frac{1}{2}&0&0\cr
                      0&0&0\cr}
                      \right).
$$

With these explicit matrices, we can calculate:
$$
  [S_1,S_2]=\frac{1}{8}\rho_3,\quad[S_2,S_3]=\frac{1}{8}\rho_1,
  \quad[S_3,S_1 ]=\frac{1}{8}\rho_2
$$
from which we deduce that $\gamma=\frac{1}{8}$ and $\delta=0.$
We also find
$$
  [S_1,\rho_2]+[\rho_1,S_2]=-4S_3 \quad {\rm etc.}
$$
  and thereby obtain $\alpha=0$ and $\beta=-4.$
  Hence, in the tetrahedral case, Nahm's equations are reduced via
  (\ref{e1}), (\ref{e2})   to
\begin{eqnarray}
\frac{dx}{ds}&=&2x^2+\frac{1}{8}y^2 \label{N3}\\
\frac{dy}{ds}&=&-4xy . \label{NN3}
\end{eqnarray}

When $T_i=x\rho_i+yS_i$ in the present situation, a straightforward
calculation gives the polynomial (\ref{eq:Nahmpoly}) as
\begin{equation}
\eta^3-\frac{1}{2}(48x^2+y^2)y\zeta(\zeta^4-1) ,
\label{spec3}
\end{equation}
so $(48x^2+y^2)y$ is a constant of integration for
(\ref{N3}),(\ref{NN3}). The requirement that the Nahm matrices are
antihermitian means that $x$ is real, and $y=iv$ with $v$ real. The
constant of integration is $(48x^2 - v^2)v = c$, where $c$ is
real, and the polynomial (\ref{spec3}) is the same as in
(\ref{A}), if we identify $a = -\frac{1}{2} c$.

Using this constant of integration and substituting in (\ref{NN3}) gives
$$\frac{\sqrt{3}dv}{\sqrt{v^4 + cv}}=-ds .$$
Putting $v = c^{1/3}\wp(u)^{-1}$, where the Weierstrass elliptic function
$\wp(u)$ satisfies the equation $\wp'(u)^2=4\wp^3(u)+4$, we find
$$
s=2{\sqrt{3}}c^{-1/3}u+K
$$
for some constant $K$. Thus from these substitutions and (\ref{NN3}) we
obtain the general solution to the equations
$$
x=\frac{c^{1/3}}{8\sqrt{3}}\frac{\wp'(u)}{\wp(u)},\quad
y=\frac{ic^{1/3}}{\wp(u)}.
$$

Now, the period lattice of $\wp(u)$ is (equilaterally) triangular,
with triangle edges along the imaginary axis, and
$\wp$ has double poles at the vertices of the triangles. (Near $u=0$,
$\wp(u)=u^{-2}+\dots \ $.) $\wp'=0$ at the
mid-points of the edges, and $\wp=0$ at the centres of the
triangles. Let $2\omega_1$
denote the real period. On the interval $[0,2\omega_1]$ there are zeros
of $\wp(u)$ at $u=\frac{2}{3}\omega_1$ and $u=\frac{4}{3}\omega_1$.
Between $\frac{2}{3}\omega_1$ and $\frac{4}{3}\omega_1$, $\wp(u)$ is
negative, reaching its minimum value $\wp(\omega_1)=-1$.

To fit the boundary conditions of a monopole, we require $T_i(s)$ to be
regular for $0<s<2$ and to have poles at $s=0$ and $s=2$ whose residues
define the irreducible three-dimensional representation of
${\frak{so}}(3)$. We can satisfy these conditions if we choose $K$ so that
$u=\frac{2}{3}\omega_1$ when $s=0$ and $u=\frac{4}{3}\omega_1$ when $s=2$,
that is
$$
s=\frac{1}{\omega_1}(3u-2\omega_1)
$$
with $2\omega_1=\sqrt{3}c^{1/3}$. Both $x$ and $y$ have simple poles at
$s=0$ and $s=2$, and one can verify that the residues of the Nahm matrices
there define an irreducible representation of ${\frak{so}}(3)$. Moreover,
$x$ and $y$ have no singularities for $0<s<2$. Since $y(2-s)=y(s)$ and
$x(2-s)=-x(s)$, it follows that $T_i(2-s)=T_i(s)^T$.

It is straightforward to see that another solution with these boundary
conditions is obtained by replacing $2\omega_1$ by $-2\omega_1$; the
result is a reflection of the monopole about the origin.

Yet another solution of Nahm's equations is obtained by choosing $K=0$,
with $u=\frac{2}{3}\omega_1$ at $s=2$. This gives an equivalent monopole,
with the same spectral curve. Although the symmetry $s \rightarrow 2-s$ is
no longer manifest, it can be obtained after a unitary transformation of
the matrices.

The elliptic curve featured in the solution is a very special one, and the
period $2\omega_1$ may be evaluated explicitly. We have (see e.g.
\cite{WhiWat})
$$
2\omega_1= \int_{-1}^{\infty}\frac{dt}{\sqrt{t^3+1}}=
\frac{1}{2\sqrt{\pi}}\Gamma(\frac{1} {6})\Gamma(\frac{1}{3})
$$
and consequently
\begin{theorem}
The curve $\eta^3+ia\zeta(\zeta^4-1)$ is the spectral curve of
a charge 3 monopole with tetrahedral symmetry if
$a=\pm \Gamma(\frac{1}{6})^3\Gamma(\frac{1}{3})^3/48
{\sqrt{3}}\pi^{3/2}.$
\label{th:tetra}
\end{theorem}

\section{Octahedral symmetry}
\label{sec:octa}
In the case of the octahedral group, $k=4$, and we  consider the
irreducible representation of $SU(2)$ on homogeneous cubic polynomials. In
a suitable basis, we can express the representation by
$$\rho_1=\left(
\matrix{0&3&0&0\cr
        -1&0&4&0\cr
        0&-1&0&3\cr
        0&0&-1&0}\right),\ \rho_2=i\left(
        \matrix{0&3&0&0\cr
                1&0&4&0\cr
                0&1&0&3\cr
                0&0&1&0}\right),\
                \rho_3=i\left(
                \matrix{3&0&0&0\cr
                        0&1&0&0\cr
                        0&0&-1&0\cr
                        0&0&0&-3}\right).$$
The matrices $S_i$ are found by polarizing the polynomial
$\zeta_1^8+14\zeta_1^4\zeta_0^4+\zeta_0^8$, and representing the resulting
three sextic polynomials as $4\times 4$ matrices. We find
$$S_1=\left(
\matrix{0&-6&0&-60\cr
        2&0&12&0\cr
        0&-3&0&-6\cr
        5/3&0&2&0}\right),\ S_2=i\left(
                            \matrix{0&-6&0&60\cr
                                    -2&0&12&0\cr
                                    0&3&0&-6\cr
                                    5/3&0&-2&0}\right)\
$$
and
$$
S_3=4i\left(
                                    \matrix{1&0&0&0\cr
                                            0&-3&0&0\cr
                                            0&0&3&0\cr
                                            0&0&0&-1}\right)$$
                     and from these expressions we obtain
$$[S_1,S_2]=-48\rho_3-8S_3\qquad{\rm and}\qquad
[S_1,\rho_2]+[\rho_1,S_2]=-6S_3\quad{\rm etc.}$$
from which it follows that
$$
\alpha=0,\quad \beta=-6,\quad \gamma=-48,\quad \delta=-8.
$$
Using (\ref{e1}) and (\ref{e2}), Nahm's equations then reduce to
\begin{eqnarray}
\frac{dx}{ds}&=&2x^2-48y^2 \label{N4}\\
\frac{dy}{ds}&=&-6xy-8y^2 . \label{NN4}
\end{eqnarray}
In this case,
 \begin{equation}
 \det(\eta +i(T_1+iT_2)-2iT_3\zeta
 -i(T_1-iT_2)\zeta^2)=\eta^4-960y(x+3y)(x-2y)^2(\zeta^8+14\zeta^4+1)
 \label{det}
 \end{equation}
 so that we have an integral of the equations
 \begin{equation}
 y(x+3y)(x-2y)^2=c.
 \label{int4}
 \end{equation}

 To solve the equations, put $x=ty$, then (\ref{NN4})
 and (\ref{int4}) give
 $$
 y^4(t+3)(t-2)^2=c\quad{\rm and} \quad
 \frac{d}{ds}\left(\frac{1}{y}\right)=2(3t+4)$$
 and hence
 $$\frac{d}{ds}((t+3)^{1/4}(t-2)^{1/2})=2c^{1/4}(3t+4).$$
 Now making the substitution $t=5w^2-3$, we obtain
 $$\frac{dw}{\sqrt{w(w^2-1)}}=4(5c)^{1/4}ds . $$
 This is an elliptic integral and is solved by setting $w=\wp(u)$ with
 $\wp'(u)^2=4\wp(u)^3-4\wp(u)$ to get
 \begin{equation}
 u=2(5c)^{1/4}s+K .
 \label{sol4}
 \end{equation}
 The general solution to the equation is then
 $$x=\frac{2c^{1/4}(5\wp^2(u)-3)}{5^{3/4}\wp'(u)},\quad
 y=\frac{2c^{1/4}}{5^{3/4}\wp'(u)} .$$
 \vskip .5cm
At $u=0$, $y$ vanishes and so if $K=0$, $x=-(5c)^{1/4}/u+\dots=-1/2s+\dots$
from  (\ref{sol4}). Thus $T_i=-\rho_i/2s+\dots$ and the boundary condition
at $s=0$ is satisfied since again $\rho$ is an irreducible representation.
Now, however, $x$ and $y$ acquire simple poles at each half-period, where
$\wp'(u)=0$.  Since $\wp'(u)^2=4(\wp(u)^3-\wp(u))$ this occurs where
$\wp(u)=0,-1,+1$. Using $\wp''(u)=6\wp(u)^2-2$, and (\ref{sol4}), we find
the residue  of $x$ to be $(5\wp^2-3)/10(3\wp^2-1)$ and that of $y$ to be
$1/10(3\wp^2-1)$.  Thus, if $\wp(u)=0$, the residue $R_3$ of the Nahm
matrix $T_3=x\rho_3+yS_3$ is
 $$\frac{i}{2}\left(
 \matrix{1&0&0&0\cr
          0&3&0&0\cr
          0&0&-3&0\cr
          0&0&0&-1}\right)$$
  which identifies the representation at this pole as $S^3V$, the
  irreducible one. On the other hand, if $\wp^2(u)=1$, the residue is
                   $$\frac{i}{2}\left(
                   \matrix{1&0&0&0\cr
                   0&-1&0&0\cr
                   0&0&1&0\cr
                  0&0&0&-1}\right)$$
  which are the weights of $V\oplus V$, which is reducible.
 It follows that only when
  $\wp(u)=0$  does the solution to Nahm's
  equations have a pole with the correct residue.

Now if $u$ is real, $\wp'(u)^2\ge0$, so $\wp^3(u)\ge\wp(u)$. Since
$\wp(u)\rightarrow +\infty$ as $u$ approaches $0$ or $2\omega_1$, where
$2\omega_1$ is the real period, the turning point $u=\omega_1$ is where
$\wp(u)=1$. Similarly if $u$ is imaginary, and the imaginary period is
$2\omega_3$, we must have $\wp(\omega_3)=-1$. Thus the required pole is at
$u=\omega_1+\omega_3=\omega_2$. From (\ref{sol4}) this is possible if
$c<0$ since we can then take the argument of $c^{1/4}$ to be $\pi/4$, and
as $s$ takes real values, $u$ lies on the line from $0$ to $\omega_2$. The
boundary condition at $s=2$ can then be satisfied if
\begin{equation}
  \omega_2=4(5c)^{1/4}
  \label{omega4}
\end{equation}
which determines the constant $c$. To put this in a more concrete form,
note that the substitution
$$
w=\frac{z-i}{1-iz}
$$
transforms the differential
$$
du=\frac{\wp'du}{\sqrt{4(\wp^3-\wp)}}=\frac{dw}{2\sqrt{w^3-w}}
$$
into
$$
(1+i)\frac{dz}{2\sqrt{z^4-1}}
$$
and so
$$
\omega_2=(1+i)\int_0^1\frac{dt}{\sqrt{1-t^4}} .
$$
{}From \cite{WhiWat}, we also have the formula

$$
\int_0^1\frac{dt}{\sqrt{1-t^4}}=\frac{1}{\sqrt{8\pi}}
\Gamma(\frac{1}{4})^2.
$$
Using this, together with (\ref{omega4}) and (\ref{det}), we obtain
\begin{theorem}
 The curve $\eta^4+a(\zeta^8+14\zeta^4+1)$ is the spectral
  curve   of a charge 4 monopole with octahedral symmetry if
$a=3\Gamma(\frac{1}{4})^8/64\pi^2$.
\label{th:octa}
\end{theorem}

\section{Icosahedral symmetry}
\label{sec:icosa}
  For the icosahedral group, $k=6$, and the representation is defined by
  \begin{eqnarray*}
  \rho_1&=&\left(
  \matrix{0&5&0&0&0&0\cr
          -1&0&8&0&0&0\cr
          0&-1&0&9&0&0\cr
          0&0&-1&0&8&0\cr
          0&0&0&-1&0&5\cr
          0&0&0&0&-1&0}\right),\\
           \rho_2&=&i\left(
  \matrix{0&5&0&0&0&0\cr
          1&0&8&0&0&0\cr
          0&1&0&9&0&0\cr
          0&0&1&0&8&0\cr
          0&0&0&1&0&5\cr
          0&0&0&0&1&0}\right),\\
         \rho_3&=&i\left(
        \matrix{5&0&0&0&0&0\cr
                0&3&0&0&0&0\cr
                0&0&1&0&0&0\cr
                0&0&0&-1&0&0\cr
                0&0&0&0&-3&0\cr
                0&0&0&0&0&-5}\right).
                \end{eqnarray*}
Using Maple, we obtain  the matrices $S_i$ by polarizing the polynomial
$\zeta_0\zeta_1^{11}+11\zeta_0^6\zeta_1^6-\zeta_0^{11}\zeta_1$:

      \begin{eqnarray*}
      S_1&=&\left(
      \matrix{0&-240&0&0&-40320&0\cr
              48&0&960&0&0&40320\cr
              0&-120&0&-1440&0&0\cr
              0&0&160&0&960&0\cr
              14&0&0&-120&0&-240\cr
              0&-14&0&0&48&0}\right),\\
      S_2&=&i\left(
       \matrix{0&-240&0&0&40320&0\cr
               -48&0&960&0&0&-40320\cr
               0&120&0&-1440&0&0\cr
               0&0&-160&0&960&0\cr
               14&0&0&120&0&-240\cr
               0&-14&0&0&-48&0}\right),\\
      S_3&=&2i\left(
           \matrix{48&0&0&0&0&40320\cr
                    0&-240&0&0&0&0\cr
                    0&0&480&0&0&0\cr
                    0&0&0&-480&0&0\cr
                    0&0&0&0&240&0\cr
                    14/5&0&0&0&0&-48}\right).
 \end{eqnarray*}
 {}From these expressions, we obtain the commutation relations
 $$[S_1,S_2]=-230400\rho_3 +480S_3\qquad{\rm etc.}$$
 and
 $$
 \qquad[S_1,\rho_2]+[\rho_1,S_2]=-10S_3\qquad{\rm etc.}.
 $$
{}From (\ref{e1}) and (\ref{e2}), and putting $z=480y$, Nahm's equations
reduce to
  \begin{eqnarray}
 \frac{dx}{ds}&=&2x^2-z^2\label{M6}\\
 \frac{dz}{ds}&=&-10xz+z^2\label{MM6}
 \end{eqnarray}
 and the polynomial (\ref{eq:Nahmpoly}) gives the constant of integration
 \begin{equation}
 c=336z(z-3x)^2(z+4x)^3.
 \label{const}
 \end{equation}

 To solve the equations, put
 $$w^3={\frac{4}{7}}(1-{\frac{3x}{z}})$$
 then (\ref{const}) becomes
 $$z^6w^6(1-w^3)^3={\frac{9c}{7^6}}={\kappa}^6$$
 and so
 \begin{equation}
 z=\frac{\kappa}{w\sqrt{1-w^3}}.
 \label{con}
 \end{equation}
 Substituting in (\ref{MM6}), we obtain
 $${\frac{d}{ds}}\left({\frac{w\sqrt{1-w^3}}{\kappa}}\right)={\frac{7}{3}}
 \left(1-{\frac{5w^3}{2}}\right)$$
 but, expanding the derivative,
 $${\frac{d}{ds}}\left(
w{\sqrt{1-w^3}}\right)={\frac{dw}{ds}}\frac{(1-5w^3/2)}{\sqrt{1-w^3}}$$
 and hence
 $$\frac{dw}{\sqrt{1-w^3}}=\frac{7\kappa}{3}ds.$$
 Making the substitution $w=\wp(u)$ where the Weierstrass elliptic function
$\wp(u)$ satisfies the equation $\wp'(u)^2=4\wp^3(u)-4$ we integrate
this by
\begin{equation}
u=i{\frac{7\kappa}{6}}s+K .
\label{uu}
\end{equation}

 Using these elliptic functions, we can compute $z$ from (\ref{con}):
 $$z=-\frac{2\kappa i}{\wp(u)\wp'(u)}$$
and from the definition of $w$,
 $$
 x=i\kappa\left( \frac{\wp'(u)}{6\wp(u)}+\frac{\wp(u)^2}{2\wp'(u)}\right).
 $$
At  $u=0$, $\wp(u)= 1/u^2+\dots$ , so $z$ vanishes and $x$ has a simple
pole. In fact, putting $K=0$ in (\ref{uu}),
 $$x= -\frac{7i\kappa}{12u}+\dots= -\frac{1}{2s}+\dots$$
 and so $T_i=-\rho_i/2s\dots$. Since $\rho$ is irreducible, this is the
required behaviour
 at $s=0$.

 The other poles occur at $u=u_0$, where $u_0$ is one of the two values
 where $\wp(u_0)=0$, and at the half-periods $u=\omega_i$.
 Near $u=u_0$, we have
 $$\wp(u)=\wp'(u_0)(u-u_0)+\dots$$
 and since $\wp'(u)^2=4\wp(u)^3-4$, we obtain
 $$x=\frac{i\kappa}{6(u-u_0)}+\dots=\frac{1}{7(s-s_0)}+\dots$$
 and
 $$
 z=\frac{i\kappa}{2(u-u_0)}+\dots=\frac{3}{7(s-s_0)}+\dots
 $$
using (\ref{uu}). Using this, we calculate the residues $R_i$ of the
matrices $T_i=x\rho_i+yS_i$ and find the  eigenvalues of $iR_3$ to be
$-1,-1,0,0,1,1$, thus identifying the representation as $S^2V\oplus S^2V$.
 At a half-period $\omega_i$,
 $$
 \wp(u)=\wp(\omega_i)+\wp''(\omega_i)(u-\omega_i)^2/2+\dots
 $$
but since  $\wp'(u)^2=4\wp(u)^3-4$, $\ \wp''(\omega_i)=6\wp(\omega_i)^2$
and this provides $x(s)$  with a  simple  pole  of  residue  $1/14$ and
$z$ a simple pole of residue  $-2/7$ at the corresponding value of $s$.
The eigenvalues of $iR_3$ are then $-1/2,-1/2,-1/2,1/2,1/2,1/2$, giving
the representation $V \oplus V \oplus V$.

It follows that only the poles at the period points give irreducible
representations. But any line joining two periods passes through a
half-period. We conclude that there is no solution to Nahm's equations of
this form which is smooth in the interval $(0,2)$, has poles at the
end-points, and   whose   residues   give   an   irreducible
representation. Consequently
  \begin{theorem}
  There does not exist a monopole of charge 6 with icosahedral symmetry.
 \label{th:icosa}
 \end{theorem}

\section{Cyclically symmetric scattering of monopoles}
\label{sec:scatt}
In this Section, we shall investigate the rational maps of monopoles with
cyclic symmetry, and shall discover some novel types
of geodesic monopole scattering. Recall that there is a $1-1$
correspondence between the maps and
monopoles. Also, cyclic or axial symmetry about the $x_3$-axis, if
present, is manifest.

The rational map
of a charge $k$ monopole takes the form
$$
R(z) = {p(z) \over q(z)},
$$
with $q$ monic of degree $k$ and $p$ of degree less than $k$.
Let $\omega = e^{2\pi i/k}$.
Consider the cyclic group of rotations about the
$x_3$-axis, $C_k$, generated by the transformation $z \rightarrow
\omega z$.
The monopole with rational map $R(z)$ is $C_k$ symmetric if $R(\omega z)$
differs from $R(z)$ only by a constant phase.  We get a class of
charge $k$ monopoles with $C_k$ symmetry for each irreducible
character of $C_k$.
Let us denote the $l$th such class of monopoles by $M_k^l \, (0 \leq
l < k)$.
These are the monopoles whose rational maps are of the form
$$
R(z) = {\mu z^l \over z^k - \nu}
$$
where $\mu$ and $\nu$ are complex parameters.  For these monopoles,
$R(\omega z) = \omega^l R(z)$.
$\, M_k^l$ is a $4$-dimensional totally geodesic submanifold of the
moduli space $M_k$, since it arises by imposing a symmetry on
the monopoles.  Its metric is also K\"ahler,
because the set of rational maps
$M_k^l$ is a complex submanifold of the set of all rational maps.

Since the strongly centred monopoles are geodesic in the moduli space,
we shall now restrict attention to rational maps of strongly centred,
$C_k$-symmetric monopoles. There is no essential loss of generality in
doing this. For a monopole with a rational map of type
$M_k^l$, the criterion for it to be strongly centred (\ref{eq:strongcent}),
reduces to
\begin{equation}
\mu^{k} \prod^k_{i=1} \left (\beta_i\right)^l  = 1 \label{eq:sc}
\end{equation}
where $\{\beta_i: i=1,\dots, k\}$ are the $k$ roots of $z^{k} - \nu = 0$.
Eq.~(\ref{eq:sc}) is equivalent to $\mu^k \nu^l = \pm 1$, with the lower
sign if both $k$ is even and $l$ odd, and the upper sign otherwise. The
magnitude of $\mu$ is $|\mu| = |\nu|^{-l/k}$, and there are $k$ choices
for the phase of $\mu$.  The rational maps we obtain are parametrised by
several surfaces of revolution. For given $k$ and $l$ there may be one or
more surfaces. For $l=0$, for example, there are $k$ distinct surfaces,
each with $\nu$ a good coordinate; $\mu$ is a distinct, and constant,
$k$th root of unity on each surface. If $l \not= 0$, and $k$ and $l$ have
highest common factor $h$, there are $h$ distinct surfaces. As arg $\nu$
increases by $2\pi$, arg $\mu$ decreases by $2\pi l/k$, so arg $\nu$ must
increase by $2\pi k/h$ for $\mu$ to return to its initial value. $\nu$ is
therefore a good coordinate on each surface, but the range of arg $\nu$ is
$2\pi k/h$.

For given $k$, and each $l$ in the range $0 \leq l < k$, let us
choose one of the surfaces just described, say, the one containing the
rational map in $M_k^l$ with $\nu = 1$ and $\mu = e^{\pi i/k}$ (if $k$ is
even and $l$ odd) or $\mu = 1$ (otherwise). Denote this surface by
$\Sigma_k^l$. If there is
another surface, for a particular value of $l$, then it is isomorphic
to $\Sigma_k^l$, as $\mu$ differs on it simply by a constant phase.
Let us now consider the geodesics on $\Sigma_k^l$, and the associated
$C_k$-symmetric monopole scattering.
The simplest geodesic is when $\nu$ moves along the real axis -- the
monopole then has no angular momentum.

On $\Sigma_k^0$ the rational maps are of the form
$$
R(z) = {1 \over {z^k - \nu}},
$$
where $\nu$ is an arbitrary complex number. $\Sigma_k^0$ is therefore
a submanifold of the space of inversion symmetric monopoles
$\M_k^0$. For $\nu =0$, the
rational map is that of a strongly centred axisymmetric charge $k$
monopole. If $|\nu|$
is large, there are $k$ well-separated unit charge monopoles at the
vertices of an $k$-gon in ${\bf R}^3$, with $x_1+ix_2$ a $k$th root of
$\nu$, and $x_3
= 0$.  The geodesic where $\nu$ moves along the entire real axis
corresponds to a
simultaneous scattering of $k$ unit charge monopoles in the
$(x_1,x_2)$ plane, where the incoming and outgoing trajectories are
related by a $\pi/k$ rotation.  The configuration is instantaneously
axially symmetric when $\nu = 0$.  This kind of symmetric planar
scattering of $k$ solitons has been observed in a number of models,
and can be understood in a rather general way \cite{KudPieZak,Dzi}.

On $\Sigma_k^l$, with $l \neq 0, \ \nu$ is a
non-zero complex number.  $\nu =
0$ is forbidden, as the numerator and denominator of $R(z)$ would have a
common factor $z^l$.  A simple geodesic is with $\nu$ moving along the
positive real axis, say towards $\nu =0$.  The rational map is
\begin{equation}
R(z) = {\iota \over \nu^{l/k}} \, {z^l \over {z^k - \nu}} \label{eq:ratmap}
\end{equation}
where $\iota = e^{\pi i/k}$ (if $k$ is even and $l$ is odd) or $\iota
= 1$ (otherwise).
Then the initial motion is again $k$ unit charge monopoles at the vertices
of a contracting $k$-gon in the $(x_1, x_2)$ plane.  As $\nu
\rightarrow 0$, the map approaches
$$
R(z) = {\iota \over \nu^{l/k}} \,  {1 \over z^{k-l}}
$$
which is the map of a charge $(k-l)$ axisymmetric monopole, centred at
the point $(0,0, (- l / 2k) \log \nu)$.  This is a positive distance along
the $x_3$-axis as $\nu$ is small.  Following an argument in [\cite{AtiHit},
pp.25-6], we deduce that the charge $k$ monopole has split up,
with one
cluster the charge $k-l$ monopole just described, and a further
cluster (or clusters) near the $x_3$-axis, but not so far up.  In fact,
there is just one other cluster, which is an axisymmetric charge $l$
monopole at a negative distance along the $x_3$-axis.  This is
seen by inverting the original monopole in the $(x_1,x_2)$ plane.  The
procedure described in Section~\ref{sec:invratmap}
 shows that the rational map (\ref{eq:ratmap})
transforms under inversion to
$$
R(z) = {\tilde\iota \over \nu^{(k-l)/k}} {z^{k-l} \over {z^k - \nu}}
$$
where $\iota \, \tilde\iota = 1$, because
$$
\begin{array}{cl}
\iota \, \tilde\iota \, (z^l / \nu^{l/k}) \, (z^{k-l} /
\nu^{(k-l)/k}) & = {z^k / \nu} \\
& = 1 \mod z^k-\nu .
\end{array}
$$
The inverted monopole therefore has an axisymmetric charge $l$
monopole cluster at $\allowbreak$
$(0,0, - ((k-l)/2k) \log \nu)$, as $\nu \rightarrow 0$, while the
original monopole has this axisymmetric charge $l$
cluster at $\allowbreak$ $(0,0,((k-l)/2k )\log\nu)$.

In the geodesic motion, $k$ unit charge monopoles come in, but the outgoing
configuration is of two approximately axisymmetric monopole clusters, of
charges $k-l$ and $l$, at distances $ld$ and $-(k-l)d$ along the
$x_3$-axis, with $d$ increasing uniformly.  This geodesic motion can, of
course, also be reversed.  The centre of mass of these clusters remains at
the origin.

If $k$ is even and $l=k/2$ then the rational maps, and the
geodesic monopole motion we have described, have an additional
inversion symmetry. $R(z) = z^{k/2}/(\nu^{1/2}(z^k - \nu))$ lies in
the space of inversion symmetric maps $\M_k^{k/2}$, and the
factor $\iota$ makes no essential difference. Consequently, the
outgoing clusters have the same charges and
equal speeds.  Since $\nu$ was assumed to be real, there is reflection
symmetry under
$x_2 \rightarrow -x_2$.  Together with the inversion symmetry, $x_3
\rightarrow -x_3$, we obtain an additional rotational symmetry, by
$\pi$ about the $x_1$-axis.  Hence, monopoles with rational maps of the
form  (\ref{eq:ratmap}) have $D_k$ symmetry if $k$ is even and $l = k/2$.
There is also $D_k$ symmetry  if $l=0$, for any $k$.

The surfaces $\Sigma_2^0$ and $\Sigma_2^1$  are the ``rounded cone''
and ``trumpet'' described in \cite{AtiHit}.  These surfaces are
not isomorphic, but the geodesics with $\nu$ real (on $\Sigma_2^0$)
and $\nu$ real and positive (on $\Sigma_2^1$) {\it are}
isomorphic.  Along the first, two unit  charge monopoles scatter
through $\pi/2$ in the $(x_1, x_2)$ plane, and along the second they
scatter through $\pi/2$ in the $(x_1, x_3)$ plane.  There are no
analogous isomorphisms in the higher charge cases.

The general geodesics on the surfaces $\Sigma_k^0$ and
$\Sigma_k^l \ (l \neq 0)$ are presumably analogous to  those on the cone
$\Sigma_2^0$ or trumpet $\Sigma_2^1$.  On $\Sigma_k^0$, they correspond
to $k$ unit charge monopoles scattering in the $(x_1,x_2)$ plane with
net orbital angular momentum.  On $\Sigma_k^l \ (l \neq 0)$, $k$ unit
charge monopoles again come in
with net orbital angular momentum.  If this is small, the geodesic
passes through the trumpet-like surface and two monopole clusters
with magnetic charges $l$ and $k-l$ emerge back-to-back on the
$x_3$-axis.  They also have opposite electric charges, which accounts,
physically, for angular momentum conservation.  If the initial angular
momentum is large, then the geodesic does not pass through the
trumpet, but is reflected, and there are $k$ outgoing unit charge monopoles
in the
$(x_1,x_2)$ plane.

What can we learn about the spectral curves of centred $C_k$-symmetric
monopoles from this discussion of rational maps? First, recall that
monopoles whose rational maps differ only by a phase have the same
spectral curves. We need therefore only consider the chosen surfaces
of rational maps, $\Sigma_k^l$, and their associated monopoles. Let us
also restrict attention to monopoles which are oriented to be
reflection symmetric under $x_2 \rightarrow -x_2$, which requires
$\nu$ to be real, and choose a fixed phase for $\mu$ as $\nu$ varies in
magnitude. This restricts us to $2k-1$ disjoint curves in the
surfaces $\Sigma_k^l$, $(\nu$ real in $\Sigma_k^0, \
 \nu$ positive and $\nu$ negative in $\Sigma_k^l \ (l
\neq 0))$, and these curves are geodesics.
It follows that among the centred curves in
$TP_1$ of the form (\ref{eq:spec}) with $C_k$
symmetry and oriented, there are $2k-1$ disjoint loci
of spectral curves. (We refer to a connected,
one-dimensional submanifold of spectral curves as a locus in the space of
curves in $TP_1$.) All these spectral curves will have real
coefficients because of the reflection symmetry. We have been unable
to determine, in general, for which parameter values a curve is a
spectral curve, but we can make some qualitative assertions, based on
knowledge of the asymptotic monopole configurations, the
axisymmetric configurations, and the monopoles with the symmetries of
the regular solids.  We restrict our remarks to the cases
$k=3$ and $k=4$.

For $k=3$, and  $l = 0, 1\ \hbox{or}\ 2$, there are five loci of
spectral curves of the form (\ref{eq:C3symm}), with $\beta$ real.
When $l=0$  there is $D_3$-symmetry, so
$\gamma = 0$.  The locus is asymptotic at both ends to $\alpha^3 =
27\beta^2$, with $\beta$ large and positive at one end, and $\beta$
large and negative at the other.  The axisymmetric monopole,  half-way
along the locus, has $\beta=0$ and $\alpha = \pi^2$. Presumably,
$\alpha$ is positive along the whole locus. The four
remaining loci, for $l=1$ and $l=2$, are isomorphic.  This is because
$\nu \rightarrow -\nu$ corresponds to a reflection $x_1 \rightarrow
-x_1$, and because the $l=2$ monopoles are obtained from $l=1$
monopoles by inversion $(x_3 \rightarrow -x_3)$.  Under the first
symmetry $\beta \rightarrow -\beta$, and under the second $\gamma
\rightarrow -\gamma$.  Each of the four loci is asymptotic at one end
to $\alpha^3 = 27\beta^2, \gamma = 0$, with $\beta$ either positive or
negative, and at the other to $\alpha = {1 \over 4}\pi^2 - 3b^2,\
\beta=0,\ \gamma
= 2b(b^2 + {1 \over 4}\pi^2)$, with $b$ either positive or negative.  These
latter parameters result from taking the product of the spectral curve
of a unit charge monopole at $(0,0,b)$ with the spectral curve of an
axisymmetric charge 2 monopole at $(0,0, - b/2),$ that is
$$
\begin{array}{cl}
P(\eta, \zeta)& = (\eta + 2b\zeta) (\eta^2 - 2b\eta\zeta + (b^2 +
{1 \over 4}\pi^2) \zeta^2) \\
& = \eta^3 + ({1 \over 4}\pi^2 - 3b^2)\eta \zeta^2 + 2b(b^2 +
{1 \over 4}\pi^2)\zeta^3 = 0 .
\end{array}
$$
Since the tetrahedrally symmetric charge $3$ monopole has $C_3$
symmetry about various axes, there must be a point on each of these four
loci corresponding to such a monopole (in four distinct orientations).
The calculations of Section \ref{sec:tetra} show that the four loci pass
through the four points (one on each locus) $(\alpha, \beta, \gamma) =
(0, \pm a, \pm 5\sqrt{2} a)$, where $a =
\Gamma(\frac{1}{6})^3\Gamma(\frac{1}{3})^3/48{\sqrt{3}}\pi^{3/2}$.
The picture of the
monopole scattering, corresponding to geodesic motion along one of
these loci, is as follows. Three unit charge monopoles come in at the
vertices of a (horizontal) equilateral triangle, moving towards its
centre. They then coalesce instantaneously into a tetrahedron,
with the base triangle oriented the same
way as the initial triangle, but somewhat below the initial plane.
Finally the tetrahedron breaks
up with the top vertex moving up and becoming a unit charge monopole,
and the base triangle descending and
becoming a toroidal charge two monopole.

In the case $k=4$, we have seven loci of spectral curves with $C_4$
symmetry.  Only three  of these are essentially different.  The four
corresponding to the rational maps with $l=1$ and $l=3$, and $\nu$
positive or negative, are isomorphic.  The $l=1$ and $l=3$ maps, and
hence the corresponding monopoles and spectral curves, are related by
inversion, and the sign of $\nu$ can be reversed  by a $\pi/4$
rotation.  The spectral curves along these four loci have no higher
symmetry than $C_4$ symmetry.

There are two isomorphic loci of spectral curves corresponding to the
$l=2$ maps.
Here there is inversion symmetry, and the spectral curves are
therefore $D_4$ symmetric and of the form (\ref{eq:D4symm}).
Reversing the sign of
$\nu$ again corresponds to a $\pi/4$ rotation, and $\beta$ changes
sign.  Along these two loci, we will find the spectral curves
corresponding to the octahedrally symmetric 4-monopole (in two
orientations). The locus with $\nu$ negative interpolates between the
asymptotic parameter values $\alpha = 4a^2, \beta = a^4, \gamma=2a^4$
with $a$ large (corresponding to four stars at $(1/\sqrt{2}) (\pm
\ a, \pm \ a, 0)$ and the asymptotic values
$\alpha = {\frac{1}{2}}\pi^2 - 8b^2,
\ \beta=0, \ \gamma = (4b^2 + {\frac{1}{4}}\pi^2)^2$ with $b$ large
(corresponding to
two axisymmetric charge $2$ monopole clusters on the $x_3$-axis).  Along
the locus,
$\alpha$ changes sign,
and when $\alpha=0$ the locus passes through the spectral curve of the
monopole with  octahedral symmetry, so $\gamma =
14\beta = 21\Gamma (\frac{1}{4})^8/32\pi^2$.

Finally, there is a single locus corresponding to the $l=0$ maps.
This interpolates between the asymptotic parameter values
$\alpha=4a^2,\ \beta= -a^4,\ \gamma = 2a^4$ and $\alpha = 4a^2, \ \beta
=a^4,\ \gamma = 2a^4$, with $a$ large, and passes through the values
$\alpha = 5\pi^2/2 ,\ \beta=0,\ \gamma = 9 \pi^4 / 16$,
corresponding to  the axisymmetric charge $4$ monopole.  Presumably,
$\alpha$ and $\gamma$  are positive along the entire locus.

In summary, our main result is that in the geodesic scattering of
monopoles of
charge $k$, with $C_k$ symmetry and angular momentum zero,
there are two kinds of motion.  First, there is the well-known
possibility of $k$ monopoles, of unit charge,  scattering in the $(x_1,
x_2)$ plane through an angle $\pi/k$.  Second, there is the novel
possibility of $k$ unit charge monopoles coming in as before, but
emerging as charge $l$ and charge $k-l$ axisymmetric monopoles moving
back-to-back along the $x_3$-axis.  $l$ can have any integer value in
the range $0 < l < k$.  In the special case of 3-monopole scattering,
with $l=1$ or $l=2$,
the field configuration passes through the
tetrahedrally symmetric 3-monopole. In 4-monopole scattering, with
$l=2$, the
configuration passes
through the octahedrally symmetric 4-monopole.

\bigskip
\noindent{\bf  Acknowledgements \hfill}

N.S.M. warmly thanks the Pure Mathematics Department of the University
of Adelaide and Professor Alan Carey for inviting him to visit, and
for financial support.  He also thanks the Physics Department of the
University of Tasmania at Hobart, and Professor Delbourgo, for
hospitality, and the Royal Society and the British Council for travel
grants.

M.K.M thanks the Australian Research Council for support and Michael
Singer for useful conversations.

We thank Paul Sutcliffe for help with the analysis of the
tetrahedrally symmetric monopole.

\end{document}